\documentclass[11pt]{article}
\usepackage{amsmath,amsthm,amssymb}
\usepackage{caption,graphicx,subfig}
\usepackage[margin=1in]{geometry}
\usepackage{paralist,array}
\usepackage[sort&compress,numbers]{natbib}
\usepackage{setspace}
\usepackage{booktabs}
\usepackage{tikz-cd}

\setcaptionmargin{0.25in}

\makeatletter\@addtoreset{equation}{section}\makeatother

\theoremstyle{definition}

\newcolumntype{M}[1]{>{\centering\arraybackslash}m{#1}}

\newcommand{\xv}{{\bf x}}
\newcommand{\yv}{{\bf y}}
\newcommand{\ev}{{\bf e}}
\newcommand{\fv}{{\bf f}}
\newcommand{\Fv}{{\bf F}}
\newcommand{\gv}{{\bf g}}
\newcommand{\hv}{{\bf h}}
\newcommand{\Mv}{{\bf M}}
\newcommand{\Av}{{\bf A}}
\newcommand{\Bv}{{\bf B}}
\newcommand{\Cv}{{\bf C}}
\newcommand{\Jv}{{\bf J}}
\newcommand{\Iv}{{\bf I}}

\newcommand*\samethanks[1][\value{footnote}]{\footnotemark[#1]}

\newcommand{\new}[1]{\textcolor{black}{#1}}

\begin{document}

\title{\huge{\textbf{Deep Learning of Conjugate Mappings}}} 

\author{
Jason J. Bramburger\thanks{Department of Applied Mathematics, University of Washington, Seattle, WA, 98105}\ \footnote{Email: jbrambur@uw.edu} 
\and 
Steven L. Brunton\thanks{Department of Mechanical Engineering, University of Washington, Seattle, WA, 98105}
\and
J. Nathan Kutz\samethanks[1]
}

\date{}
\maketitle

\begin{abstract}
Despite many of the most common chaotic dynamical systems being continuous in time, it is through discrete time mappings that much of the understanding of chaos is formed. Henri Poincar\'e first made this connection by tracking consecutive iterations of the continuous flow with a lower-dimensional, transverse subspace. The mapping that iterates the dynamics through consecutive intersections of the flow with the subspace is now referred to as a Poincar\'e map, and it is the primary method available for interpreting and classifying chaotic dynamics. Unfortunately, in all but the simplest systems, an explicit form for such a mapping remains outstanding.  
This work proposes a method for obtaining explicit Poincar\'e mappings by using deep learning to construct an invertible coordinate transformation into a conjugate representation where the dynamics are governed by a relatively simple chaotic mapping. 
The invertible change of variable is based on an autoencoder, which allows for dimensionality reduction, and has the advantage of classifying chaotic systems using the equivalence relation of topological conjugacies.  Indeed, the enforcement of topological conjugacies is the critical neural network regularization for learning the coordinate and dynamics pairing. We provide expository applications of the method to low-dimensional systems such as the R\"ossler and Lorenz systems, while also demonstrating the utility of the method on infinite-dimensional systems, such as the Kuramoto--Sivashinsky equation.
\end{abstract}

\section{Introduction} 

Much of our modern understanding of chaotic dynamical systems comes from tracking the intersection of trajectories with a lower-dimensional subspace that is transverse to the flow. Such a subspace is typically called a Poincar\'e section in honour of Henri Poincar\'e who pioneered this method in the late nineteenth century. Chaotic systems can be characterized by considering the iterates between successive intersections of a single trajectory with the Poincar\'e section. With such a mapping, now termed a Poincar\'e or first-return mapping, it is possible to forecast chaotic systems and also to better understand the structure of a chaotic attractor. Unfortunately, there is almost no system of interest for which we have an explicit representation of such a mapping. Thus, studies of chaotic systems often involve analyzing simplified heuristic mappings that, in most cases, cannot be directly related back to the systems that motivated their study. Prominent examples of this include the relationship between the R\"ossler system and the logistic map~\cite{Rossler}, and the Lorenz system and the H\'enon map~\cite{Henon}. The goal of this work is to overcome this disconnect by using deep learning to find an invertible coordinate transformation that allows for the explicit representation of the chaotic dynamics on the Poincar\'e section while retaining a direct connection back to the original system.

The power of modern computing and the increasing availability of time-series data has resulted in a number of data-driven discovery methods for identifying nonlinear dynamical systems models~\cite{Schmidt2009science,SINDy,CINDy,Champion,PINNs1,PINNs2,Rudy,Schaeffer,Yao,Brunton2021koopman}. The {\em sparse identification of nonlinear dynamics} (SINDy) method~\cite{SINDy} employs sparse regression to best fit the time derivative of time-series data to a sparse linear combination of candidate functions. SINDy has been adapted to discover partial differential equations~\cite{Rudy,Schaeffer}, multiscale data~\cite{BramburgerSlow,Champion2}, boundary value problems~\cite{SINDyBVP}, conservation laws~\cite{Kaiser}, implicit dynamics~\cite{SINDyPI}, conservative fluid dynamics models~\cite{Loiseau2017jfm}, stochastic systems~\cite{boninsegna2018sparse,callaham2020nonlinear}, and, most important for our work here, Poincar\'e maps~\cite{Bramburger}. In the case of Poincar\'e maps, it was shown that SINDy can be used to discover an accurate mapping when the attractor is a fixed point, periodic, or even quasiperiodic, but the method fails in many cases to provide models that govern chaotic section data. These discovered Poincar\'e maps may be inaccurate if the library of candidate functions contains insufficient terms to describe their dynamics, or when the section data is high-dimensional. The former can potentially be remedied by continuing to add terms to the library, but at present this amounts to a heuristic search that may result in overfitting and ill-conditioned numerical implementation. 

To remedy the issue of discovering mappings from chaotic data, deep learning will be leveraged for the joint discovery of a simplifying coordinate system and the simplified dynamical systems model in these coordinates~\cite{lusch2018deep,Champion,gin2019deep}.  Deep learning provides a flexible architecture for the representation of data, which has led to its significant integration into the physical and engineering sciences~\cite{goodfellow2016deep,pathak2018model,battaglia2018relational,Wehmeyer2018jcp,Mardt2018natcomm,lu2019deepxde,bar2019learning,cranmer2019learning,Raissi2019jcp,brunton2019data,Duraisamy2019arfm,Noe2019science,Brunton2020arfm,li2020fourier,raissi2020science,cranmer2020lagrangian,lee2020model,li2020multipole,li2020neural,rackauckas2020universal,kochkov2021machine,Yannis,Yannis2,Yannis3}.
Specifically, a feedforward, autoencoder neural network structure is used to implement the concept of topological conjugacy, wherein a nonlinear mapping
\begin{equation}
    \xv_{n+1} = \fv(\xv_n)
\end{equation}
is shown to generate equivalent dynamics to another mapping
\begin{equation}
    \yv_{n+1} = \gv(\yv_n) 
\end{equation}
if and only if there exists a homeomorphism $\hv$ such that
\begin{equation}\label{Conj}
    \fv = \hv^{-1}\circ\gv\circ\hv.   
\end{equation}
The relationship \eqref{Conj} gives that the relevant topology of the phase spaces of $\fv$ and $\gv$ is preserved under the homeomorphism $\hv$. In particular, fixed and periodic orbits of $\fv$ and $\gv$ are mapped into each other in a one-to-one fashion, while their stability is also preserved under the action of $\hv$. Topological conjugacies are often employed to establish local results, such as the Hartman--Grobman theorem~\cite{Guckenheimer}, to better understand the dynamics in a small region of interest in the phase space. The difference here is that we seek a {\em global} conjugacy that maps the entire phase space of $\fv$ into the entire phase space of $\gv$. Therefore, the topological organization of trajectories in phase space is the same before and after applying the homeomorphism $\hv$.   

The \new{neural} network exploits the concept of a topological conjugacy by circumventing discovering $\fv$ directly and instead discovers this invertible change of variable $\hv$ along with the `simple' mapping $\gv$ that governs iterates of the transformed observations. Here the term simple can be defined by the user and is informally determined by which and how many candidate functions are used to build the mapping $\gv$. Parsimonious representations are often desired due to interpretability and generalizability. Thus, unlike~\cite{lusch2018deep,gin2019deep} where a coordinate system is learned to pair with linear (Koopman) dynamics, or~\cite{Champion} where a coordinate system is learned to pair with a parsimonious (SINDy) representation of the dynamics, here the coordinate system is learned that pairs with a conjugacy mapping. \new{Therefore, although there are computational methods for determining homeomorphisms between two given mappings~\cite{Bollt,Bollt2}, here we take a data-driven approach and discover the homeomorphism and the conjugate mapping simultaneously.} The advantage the network presents is that the candidate functions can be kept relatively basic, while the network performs the difficult task of transforming the data to be fit by the simple mapping. Discovering simple mappings $\gv$ that generate equivalent dynamics to the unknown map $\fv$, is a data-driven generalization of what has been done for decades with chaotic dynamics: providing a simple mapping that can be used to understand the chaotic dynamics of a system of interest. The difference now is that through the discovered homeomorphism $\hv$, a direct correspondence is learned between the heuristic mapping and the original system dynamics.  

Employing neural networks to discover conjugate mappings provides the following benefits: 
\begin{enumerate}
	\item Improved forecasting of chaotic dynamics. 
	\item Dimensionality reduction.
	\item Actionable and interpretable mappings.
\end{enumerate}
Point (1) is achieved through the homeomorphism $\hv$ by transforming into a coordinate system where the section data is better fit by a sparse combination of the SINDy library; without $\hv$, the original data often doesn't admit a simple representation in a standard library. 
An improved fit to $\gv$ in the transformed variables results in an improved fit for $\fv = \hv^{-1}\circ\gv\circ\hv$ in the original coordinates. 

Point (2) comes from the fact that the homeomorphism $\hv$ only operates on section data intersected with the chaotic attractor. Therefore, it may be the case that observations are made up of a low-dimensional chaotic attractor embedded in a high-dimensional phase space. The neural network can obtain a change of coordinates that removes redundancies in the observation data and describe the chaotic attractor with the fewest degrees of freedom possible. To estimate the dimension of the attractor we will use the Kaplan--Yorke dimension~\cite{Kaplan} which is derived from the Lyapunov spectrum of the system. This is motived in part by Kaplan and Yorke's conjecture that for many chaotic systems their measure of dimension is equal to the information dimension of the attractor~\cite{Kaplan2}, and that it is useful in approximating an attractor's fractal and Hausdorff dimensions~\cite{Kuznetsov}.  

Point (3) is a consequence of the limited library used to build the conjugate mapping $\gv$. The main advantage of having an explicit form for the Poincar\'e mapping is the ability to identify its periodic/recurrent orbits. It is known that infinitely many unstable periodic orbits (UPOs) densely fill a chaotic attractor, and these UPOs manifest themselves as recurrent points in the Poincar\'e map. In this way, UPOs form the skeleton of the attractor and a chaotic transient can be visualized as closely following one UPO and then randomly switching to follow another, {\em ad infinitum}. This approach to understanding chaotic dynamics is further emphasized by noting that concepts such as Lyapunov exponents, entropy, and fractal dimensions can be expressed in terms of weighted sums over the embedded UPOs~\cite{UPO1,UPO2,UPO3}. Taking advantage of the equivalence of the dynamics generated by the unknown Poincar\'e map $\fv$ and the map $\gv$, one may identify UPOs in $\gv$ and then map them bijectively to UPOs in $\fv$ using $\hv^{-1}$. This results in a quick and easy method for identifying the location of UPOs in a Poincar\'e section that is potentially much simpler than other data-driven methods~\cite{SoUPO}. Furthermore, the identified UPOs in the section can be used to seed shooting methods to find the corresponding UPOs in the full continuous-time dynamical system, after which a fast fixed point iteration will converge to the true UPO. 

The combination of dimensionality reduction and obtaining explicit mappings holds great promise for exploiting UPOs in the flow of high- and infinite-dimensional systems, such as partial differential equations (PDEs). First, by identifying UPOs in the Poincar\'e section, the \new{Ott-Grebogi-Yorke} (OGY) method can be implemented for controlling the output of a chaotic system~\cite{OGY}. This method relies on applying small, precise parameter kicks each time the flow crosses the section to keep trajectories close to a chosen UPO. In this way, one no longer has a chaotic output from the system, but one that is periodic and therefore predictable. Second, much like in low-dimensional chaotic systems, UPOs have been shown to guide the spatio-temporal chaos observed in turbulent fluid models~\cite{Fluid1,Fluid2,Fluid3,Fluid4,Fluid5,Fluid6,Fluid7,Fluid8,Fluid9}. The work of Yalniz et al. exploits the transient visits of turbulent solutions to the 3D Navier--Stokes equations to neighbourhoods of UPOs by reducing the system to a Markov chain where each node corresponds to the neighbourhood of a UPO~\cite{Yalniz}. Since finding UPOs is a major technical challenge for implementing this method on a wide variety of fluid flows, we are hopeful that the methods presented here will be used to facilitate this analysis. To this end, the utility of the methods is demonstrated for chaotic infinite-dimensional flows on the Kuramoto--Sivashinsky PDE, a prototypical chaotic PDE used to understand turbulent flow. 

This paper is organized as follows. In Section~\ref{sec:Framework} the relevant mathematical concepts that form the basis for the deep conjugacy computational algorithms are presented. This includes defining Poincar\'e sections, conjugacies, and Lyapunov exponents, while also presenting numerical procedures for seeding initial guesses for UPOs and controlling chaos. Section~\ref{sec:Architecture} contains a description of the network architecture which is entirely designed to mimic the conjugacy $\fv = \hv^{-1}\circ\gv\circ\hv$, as described above. Various aspects of the method are illustrated in Section~\ref{sec:Low} on three distinct three-dimensional chaotic systems. Each system is strategically chosen to demonstrate that the neural network can be used to provide nonlinear changes of variables to better fit polynomial mappings, target certain types of mappings to confirm long-standing conjectures, and reduce the dimension of the observed variables down to the dimension of the attractor. Once the various advantages to the neural network are illustrated, two infinite-dimensional systems are considered in Section~\ref{sec:Infinite}. These systems are the Kuramoto--Sivashinsky PDE and the Mackey--Glass delay-differential equation. In both cases we show that the Poincar\'e section dynamics can be accurately captured by simple quadratic mappings, resulting in the extraction of a number of UPOs in the case of the former and numerical evidence for a decades old conjecture in the case of the latter. We conclude in Section~\ref{sec:Discussion} with a discussion of our findings and some avenues for future work.

\section{Mathematical Framework}\label{sec:Framework} 

The construction of deep conjugacy relations requires a diversity of mathematical methods.  As such, the following section highlights the critical mathematical methods leveraged for the task, allowing for an implementation and interpretation of the neural network architecture.  We compliment this theoretical presentation with brief discussions of numerical implementations.    

\subsection{Poincar\'e Maps}

Many of our modern methods for understanding periodic, recurrent, and chaotic dynamics comes from the formative work of Henri Poincar\'e in the late nineteenth century. He proposed that to better analyze the flow of a continuous time dynamical system
\begin{equation}
	\dot \xv = \Fv(\xv), \quad \Fv:\mathbb{R}^D \to \mathbb{R}^D
\end{equation}  
one could define a transverse-to-the-flow hypersurface $X\subset \mathbb{R}^D$ and simply track the values of the trajectories $\xv(t)$ each time they intersect the surface. This takes the dynamics of a continuous system with a $D$-dimensional phase space to a discrete collection of iterates in the (at most) $(D-1)$-dimensional hypersurface $X$. For a trajectory $\xv(t)$ we obtain a sequence
\begin{equation}
	\{\xv_n := \xv(t_n)|\ \xv(t_n)\in X,\ t_n > t_{n-1}\}
\end{equation}  
consisting of the successive intersections of the trajectory with the surface $X$. One may then consider an iterative scheme 
\begin{equation}\label{Discrete}
    \xv_{n+1} = \fv(\xv_n).
\end{equation}
\new{The above mapping $\fv$} is referred to as a {\em Poincar\'e map} or first return map, with the hypersurface $X$ typically referred to as a {\em Poincar\'e section}. 

One major benefit in moving from a continuous time dynamical system to a Poincar\'e map is the reduction in the phase space since we may discard information about the trajectories away from their intersection with the lower-dimensional Poincar\'e section. Beyond this, periodic orbits of a continuous time system manifest themselves as a recurrent set of points in the Poincar\'e section, and therefore Poincar\'e maps provide an accessible method of identifying periodic orbits in a continuous flow field. This approach is predicated on the idea that one has access to a closed form Poincar\'e mapping $\fv:X\to X$, which traditionally has been a difficult problem. Recent methods have leveraged data-driven model discovery methods to obtain such mappings~\cite{Bramburger}, although these methods become less accurate as the dimensionality $D\geq1$ of the phase space becomes large. In what follows we will outline how we may leverage state-of-the-art machine learning techniques to perform dimensionality reduction and nonlinear coordinate transforms to obtain simple mappings which generate equivalent dynamics to an unknown Poincar\'e map.

\subsection{Conjugacies and Dynamical Equivalences}

The previous subsection motivates our study of discrete dynamical systems of the form \eqref{Discrete}. In practice, the dynamics generated by the function $\fv: X \to X$ may be difficult to analyze, and so it is often strategic to establish an equivalence of the dynamics with another mapping $\gv: Y \to Y$ whose dynamics are better understood. Such an equivalence is achieved by obtaining a homeomorphism, or {\em conjugacy}, $\hv:X \to Y$ so that $\gv \circ \hv = \hv\circ \fv$. Importantly, this gives that iterating \eqref{Discrete} is equivalent to applying $\hv$ to $\xv_n$, then $\gv$, and then $\hv^{-1}$, or $\fv = \hv^{-1}\circ \gv\circ \hv$. The situation is summarized with the commutative diagram:
\begin{equation}\label{ComDiag}
	\begin{tikzcd}
X \arrow{r}{\fv} \arrow[swap]{d}{\hv} &  X\arrow{d}{\hv} \\
Y\arrow{u} \arrow{r}{\gv} & Y\arrow{u}
\end{tikzcd}
\end{equation} 
Hence, iterates governed by $\fv$ lie in one-to-one correspondence  with iterates of $\gv$, mapped into each other by $\hv$ and $\hv^{-1}$. If such an $\hv$ can be found, we say that $\fv$ and $\gv$ are {\em topologically conjugate}, or simply {\em conjugate}. Conjugacies are used to demonstrate dynamical equivalence of two mappings and they form an equivalence relation of the set of all continuous surjections of a topological space to itself. We note that conjugacies have also been referred to as factor maps, particularly in the Koopman operator literature \cite{Mezic,Mezic2}. 

In practice finding a conjugacy between mappings is exceedingly difficult, but potentially tremendously rewarding. For example, suppose that $\fv$ is a Poincar\'e mapping from the intersection of the chaotic attractor with the Poincar\'e section back into itself. It may be the case that the dynamics are difficult to analyze, and so establishing a conjugacy of $\fv$ with a simpler mapping $\gv$ would help to better understand the chaotic dynamics of $\fv$. Importantly, periodic orbits of $\fv$ and $\gv$ lie in one-to-one correspondence and their stability properties are transferred between the mappings using the conjugacy. As we will see in what follows, we can identify periodic orbits in the (purportedly) simpler mapping $\gv$, use $\hv^{-1}$ to transfer them over to points in the Poincar\'e section as periodic orbits of $\fv$, and use these points to either seed initial guesses for identifying UPOs or implement an algorithm to tame the output of the chaotic system. This method of identifying UPOs and the algorithm for controlling chaos are outlined in the following two subsections.

\subsection{Discovery and Construction of UPOs}\label{subsec:UPO}

As discussed in the introduction, UPOs can be thought of as forming the skeleton of a chaotic attractor. Therefore, to best understand the geometry and dynamics of the attractor, we seek to identify these periodic orbits. To begin, a $T$-periodic solution of the ordinary differential equation
\begin{equation}\label{ContinuousUPO}
	\dot \xv = \Fv(\xv)	
\end{equation}
satisfies $\xv(t+T) = \xv(t)$ for all $t \in \mathbb{R}$. In almost every application the period $T > 0$ is unknown, and therefore we will introduce $t = sT$ so that the unknown $T$ becomes a parameter in the system. Then, $T$-periodic solutions of \eqref{ContinuousUPO} are equivalent to $1$-periodic solutions of 
\begin{equation}\label{ContinuousUPO2}
	\frac{\mathrm{d} \xv}{\mathrm{d}s} = T\Fv(\xv).
\end{equation}
Following~\cite{PeriodicOrbits}, the search for $1$-periodic UPOs in \eqref{ContinuousUPO2} is equivalent to seeking a solution $(\xv,T)$ to the system of (nonlinear) equations
\begin{equation}\label{UPOFinder}
	\begin{split}
		\frac{\mathrm{d} \xv}{\mathrm{d}s} - T\Fv(\xv) &= 0 \\
		\xv(1) - \xv(0) &= 0\\
		\int_0^1 \langle \dot{\xv}_0(s),\xv_0(s) - \xv(s) \rangle \mathrm{d}s &=0
	\end{split}
\end{equation} 
where ${\xv}_0(s)$ comes from an initial guess $(\xv_0,T_0)$ at the solution $(\xv,T)$. The first two equations guarantee that the solution $\xv(s) = \xv(t/T)$ is a $1$-periodic solution of \eqref{ContinuousUPO2}, and therefore a $T$-periodic solution of \eqref{ContinuousUPO}. Since the system is autonomous, it follows that any time shift of a solution is also a solution, leading to a non-uniqueness of solutions to the first two equations in \eqref{UPOFinder}. To remedy this, the third equation fixes a specific time-shift, based on the initial guess, and in turn makes the solution unique.

The nonlinear equations above can be discretized in $s$ to produce a nonlinear system of equations that can be solved using any number of root-finding algorithms. We elect to use the built-in MATLAB routine fsolve. The major difficulty of searching for UPOs then becomes finding an appropriate initial condition that converges to a desired UPO. Recall that a chaotic attractor has infinitely many UPOs that densely fill the attractor, and therefore \eqref{UPOFinder} may have infinitely many solutions. In the previous subsection we described a method for obtaining the location of periodic points in the Poincar\'e section using a conjugate mapping. These periodic points correspond to the UPOs intersecting the section, and therefore we may use these periodic points to initialize our root-finding routine. Indeed, given a sequence of periodic points $\{\bar\xv_1,\bar\xv_2,\dots,\bar\xv_N\}$ in the Poincar\'e section, we may consider the initial value problems
\begin{equation}
	\begin{split}
		&\dot \xv = \Fv(\xv) \\ 
		&\xv(0) = \bar\xv_k,\ k = 1,\dots,N.
	\end{split}	
\end{equation} 
We flow each of these trajectories until they reach the Poincar\'e section again. This produces a sequence of trajectories starting from an iterate of the UPO in the Poincar\'e section that flow forward until reaching the section again, likely reaching a point near the next iterate in the section. Sequentially concatenating all of these trajectories together produces an approximation of the UPO, although with jump discontinuities each time the trajectories reaches the section due to the concatenation. The length of this initial guess further produces an initial guess for the period, $T_0$. 

\new{In our present formulation we seek to employ a neural network not to discover the Poincar\'e mapping explicitly, but a mapping that is conjugate to it along with the homeomorphism that transfers one between the mappings. Therefore, we will see that with an explicit conjugate mapping, $\gv$, we may extract a sequence of its periodic points $\{\bar\yv_1,\bar\yv_2,\dots,\bar\yv_N\}$ and then use the inverse homeomorphism $\hv^{-1}$ to obtain the desired sequence $\{\bar\xv_1,\bar\xv_2,\dots,\bar\xv_N\}$ such that $\bar\xv_k = \hv^{-1}(\yv_k)$ for all $k = 1,\dots N$. We will demonstrate with our examples below that this process produces a fast and accurate method for obtaining UPOs in chaotic systems.}

\subsection{Controlling Chaos}\label{subsec:Control}

In their seminal work~\cite{OGY}, Ott, Grebogi, and Yorke described a method for exploiting the UPOs that densely fill a chaotic attractor to tame the output of a parameter-dependent system. This approach to {\em controlling chaos} relies on finding small parameter manipulations to keep the flow of a chaotic system near a UPO, thus making the output periodic and therefore predictable. We briefly illustrate by assuming that 
\begin{equation}\label{PMapOGY}
	\xv_{n+1} = \fv(\xv_n,\mu)
\end{equation}
is a parameter-dependent Poincar\'e mapping of a chaotic dynamical system. Suppose $\{\bar\xv_1,\bar\xv_2,\dots,\bar\xv_N\}$ is a periodic sequence of the Poincar\'e mapping at a parameter value $\mu = \bar\mu$. Linearizing system \eqref{PMapOGY} about each of the $N$ points $(\bar \xv_k,\bar\mu)$ means that locally we have 
\begin{equation}
	\xv_{n+1} \approx \bar\xv_{k+1} + \Av_k(\xv_n - \bar\xv) + \Bv_k(\mu_n - \bar\mu),
\end{equation}
where $\xv_{N+1} = \xv_1$ by periodicity. Here the matrices $\Av_k$ and $\Bv_k$ are the Jacobian matrices of the Poincar\'e map evaluated at each pair $(\bar\xv_k,\bar\mu)$. Since the periodic orbit is assumed to be unstable, the matrix $\Av_1\Av_2\cdots\Av_N$ has at least one eigenvalue outside the unit circle in the complex plane. Following~\cite{Romeiras}, we can employ the pole-placement method to stabilize the periodic orbit. This is achieved by introducing small parameter manipulations at each iterate of the form
\begin{equation}
	\mu_n = \begin{cases}
		\bar\mu + \Cv_k(\xv_n - \bar\xv_k) & |\xv_n - \bar\xv_k| \leq \eta \\
		\bar\mu & \mathrm{otherwise}
	\end{cases}
\end{equation}
where $\eta> 0$ is small. Notice that with $\mu_n$ as above, near each $\bar\xv_k$ we have   
\begin{equation}
	\bar\xv_{n+1} \approx \bar\xv_{k+1} + (\Av_k + \Cv_k\Bv_k)(\xv_n - \bar\xv_k),	
\end{equation}
and so the matrices $\Cv_k$ are chosen in such a way to ensure that all eigenvalues of
\begin{equation}
	(\Av_1 + \Cv_1\Bv_1)(\Av_2 + \Cv_2\Bv_2)\cdots(\Av_N + \Cv_N\Bv_N)
\end{equation} 
lie inside the unit circle of the complex plane. With such $\Cv_k$ the UPO is now stable. Note that the choice of $\Cv_k$ can be automated using linear matrix inequalities~\cite{BramburgerUPO}. 

Although the above controlling chaos algorithm works perfectly on paper, it suffers from the fact that it requires an explicit description of the Poincar\'e map \eqref{PMapOGY}. In fact, much work has been done to circumvent this issue and it is understood that one may implement such a control algorithm with only the locations of the UPOs in the section. \new{Indeed, given a sequence $\{\bar\xv_1,\bar\xv_2,\dots,\bar\xv_N\}$ corresponding to a UPO, the derivatives in $\Av_k$ and $\Bv_k$ can be estimated directly from the collected section data. For example, if the governing equations are known, one may estimate derivatives evaluated at a pair $(\bar\xv_k,\bar\mu)$ in the following way. Initialize a trajectory of the system with initial condition $(\xv_k \pm h\ev_j,\bar\mu)$, where $\ev_j$ is the $j$th coordinate vector and $h > 0$ is small. We flow the trajectory forward until reaching the section again, so that this terminal value gives the value of $\fv(\xv \pm h\ev_j,\bar\mu)$. Hence, we may estimate 
\begin{equation}
	\partial_j \fv(\bar\xv_k,\bar\mu) \approx \frac{1}{2h}[\fv(\bar\xv_k+ h\ev_j,\bar\mu) - \fv(\bar\xv_k-h\ev_j,\bar\mu)].
\end{equation}
A similar procedure can be used to estimate the derivative in $\mu$ with the initial conditions $(\xv,\bar\mu \pm h\ev_j)$}.

\new{Therefore, the main challenge to applying the control of chaos algorithm to a broad range of physical systems is identifying the location of the UPOs in the Poincar\'e section. Although there are methods to extract these UPOs from section data~\cite{SoUPO}, UPOs can be identified quickly and efficiently by applying the inverse homeomorphism $\hv^{-1}$ to the UPOs of the discovered conjugate mapping. Since the goal of this work is to identify an explicit, relatively simple, and low-dimensional conjugate mapping, it follows that extracting its UPOs should be able to be automated with relative ease. In Section~\ref{sec:Low} we discuss this efficiency as it applies to both the R\"ossler and Gissinger systems, while also providing MATLAB implementations for controlling these systems in the code repository associated to the manuscript.}

\subsection{Lyapunov Exponents}\label{subsec:LEs} 

A standard approach to quantify chaotic behaviour is by analyzing the Lyapunov exponents of an attractor, which we now describe. We start with an ordinary differential equation 
\begin{equation}\label{Continuous}
	\dot \xv = \Fv(\xv)
\end{equation}
and a solution to the system $\xv(t)$. The Lyapunov exponents (LEs) describe the rate of contraction or separation of trajectories that start nearby $\xv(t)$. Precisely, let $\Jv(t) := \mathrm{D}\Fv(\xv(t))$ be the Jacobian matrix of $\Fv$ at $\xv(t)$ for each $t$. Then, we can define the evolution of the tangent vectors $\Mv(t)$ to be the trajectory via the differential equation equation
\begin{equation}
	\dot \Mv = \Jv(t)\Mv
\end{equation}  
with the initial condition $\Mv(0) = \Iv$, the identity matrix. The LEs are defined to be the eigenvalues of the matrix 
\begin{equation}
	\Lambda := \lim_{t\to\infty}\frac{1}{2t}\log(\Mv(t)\Mv(t)^T).
\end{equation} 
The limiting matrix, $\Lambda$, is guaranteed to exist by Oseledets theorem and is a symmetric matrix, making the eigenvalues real. When $\xv(t)$ is a steady-state solution, i.e. independent of $t$, then the LEs are just the real part of the eigenvalues of the constant matrix $\Jv$. Similarly, when $\xv(t)$ is periodic, then the LEs correspond to the real part of the Floquet exponents. 

Although the LEs can be used to understand the stability of steady-state and periodic solutions, they are primarily employed to quantify systems with chaotic behaviour. Throughout we can order the real-valued LEs such that 
\begin{equation}
	\lambda_1 \geq \lambda_2 \geq \dots \geq \lambda_D,
\end{equation}
where $D \geq 1$ is the dimension of the system. For chaotic systems, Oseledets theorem guarantees that the LEs are independent of $\xv(0)$, and chaos is typically indicated by the presence of at least one positive LE. Furthermore, one LE is always zero, corresponding to the eigenvector in the direction of the flow of $\xv(t)$. Throughout we will calculate the LEs using the work of~\cite{Wolf}, for which the LEs will only be reported to two decimal places due to potential numerical inaccuracy.  

The significance of the LEs is that they can be used to estimate the dimension of a chaotic attractor. Indeed, the Kaplan--Yorke dimension (sometimes Lyapunov dimension) is defined as 
\begin{equation}
	D_{KY} = k + \sum_{n = 1}^k \frac{\lambda_i}{|\lambda_{k+1}|},
\end{equation}
where $k$ is the maximum integer such that the sum of the $k$ largest LEs is nonnegative. The Kaplan--Yorke dimension is conjectured to be the dimension of the chaotic attractor and from this conjecture we will use it as an approximation of the dimension of the latent variables inside the network. This implies that although we may have a description of the iterates inside a Poincar\'e section using a large number of variables, it could be the case that the attractor dimension is small, as measured by the Kaplan--Yorke dimension. Thus, a conjugacy between the data collected in the section and a new map $\gv$ does not need to use the same number of variables as the collected data since the conjugacy is only established from attractor to attractor. This will allow for a notion of dimensionality reduction from the observed variables to the latent variables governed by the map $\gv$. The dimension of these latent variables, $d\geq 1$, will be given by
\begin{equation}
	d \approx D_{KY} - 1,
\end{equation}
where some rounding will be necessary since $D_{KY}$ is rarely an integer. Recall that $D_{KY}$ is the dimension of the attractor in the continuous system \eqref{Continuous}, and so we subtract one dimension as we are moving to discrete iterates in the Poincar\'e section.

\section{Network Architecture}\label{sec:Architecture} 

\begin{figure} 
\center
\includegraphics[width = 0.75\textwidth]{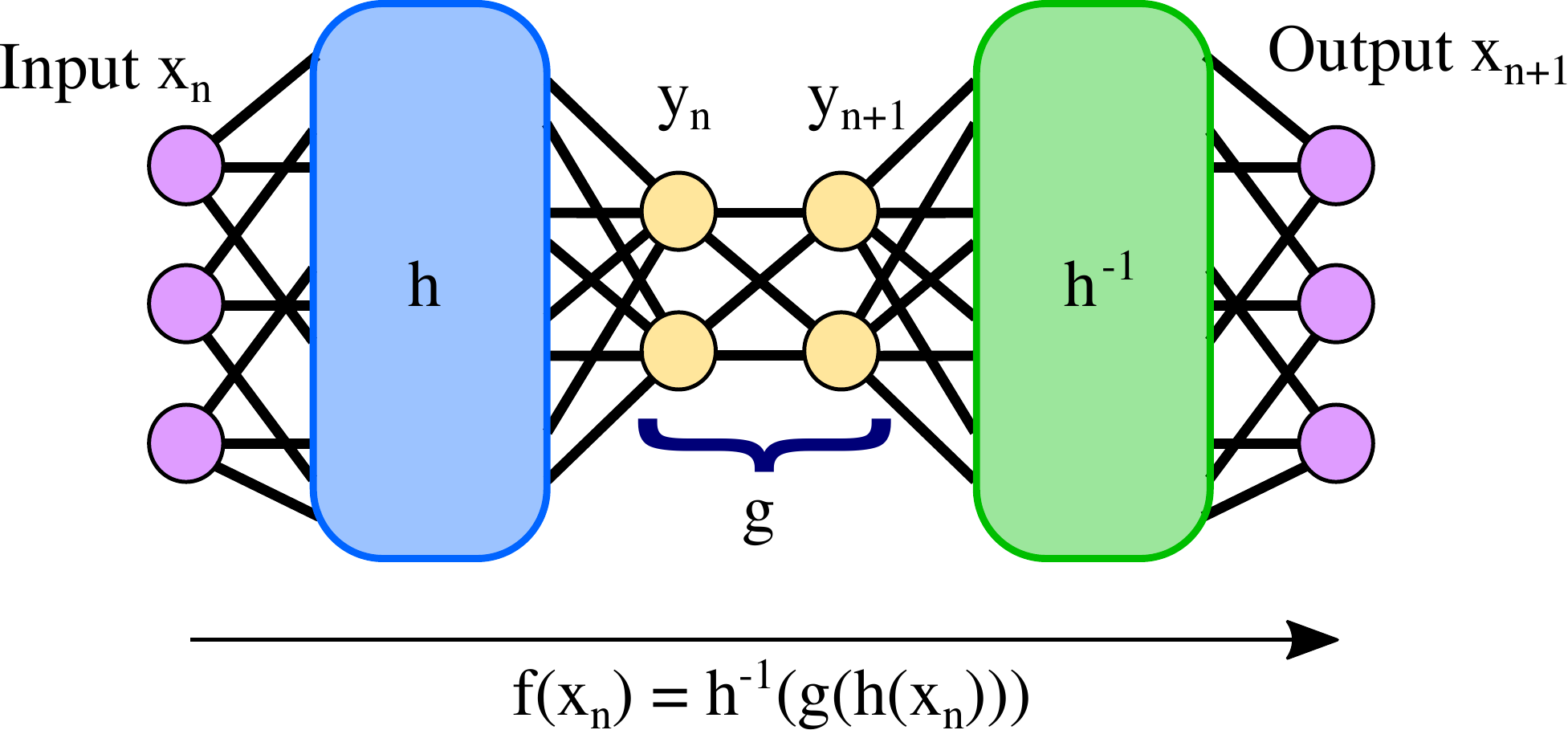} 
\caption{Diagram of the deep learning architecture for discovering nonlinear mappings $\fv$. The network first conjugates the inputs $\xv_n$ into the latent variables $\yv_n = \hv(\xv_n)$. Then, the latent variables are stepped forward using a discovered mapping $\gv:\mathbb{R}^d\to\mathbb{R}^d$. Finally, the outputs of the unknown mapping $\fv$ are recovered by applying $\hv^{-1}$ to the stepped latent variables $\gv(\hv(\xv_n)) = \gv(\yv_n)$.}
\label{fig:Network}
\end{figure} 

The proposed network architecture is based on the concept of obtaining a conjugate mapping and is illustrated in Figure~\ref{fig:Network}. This neural network takes the form of a feedforward autoencoder. The general idea is that if $\{\xv_n\} \subset \mathbb{R}^{D-1}$ is a set of observations obtained from intersections of the flow with the Poincar\'e section, then instead of attempting to identify a mapping
\begin{equation}
	\xv_{n+1} = \fv(\xv_n),
\end{equation}	
we use the the network to find a homeomorphism $\hv$ and a conjugate mapping $\gv$ such that $\fv = \hv^{-1}\circ \gv\circ \hv$. The neural network combines both nonlinear changes of coordinates and dimensionality reduction to arrive at a latent mapping $\gv$ as a function of the latent variables $\yv\in\mathbb{R}^d$. Like other data-driven discovery algorithms, the function $\gv$ is expanded as a basis of library functions. The coefficients in front of each basis function are set as variables in the network which are tuned along with the weights of the network. The conjugacy function mimics an autoencoder with $\hv$ and $\hv^{-1}$ trained as an encoder and decoder, respectively. We may collect training data from one long trajectory of the nonlinear dynamical system since the transitivity of orbits on a chaotic attractor can be relied upon to densely fill the attractor. This means that the longer the trajectory, the more we have mapped out the entire attractor in the observable variables in the section, thus densely filling out the domain of the conjugacy function $\hv$. 

The goal is that the conjugate/latent mapping $\gv:\mathbb{R}^d\to\mathbb{R}^d$ is simple, as defined by the user. For example, polynomials are some of the simplest basis functions that a function can be expanded in, but it is unlikely that the unknown function $\fv$ is itself a polynomial. With this in mind, we can create a library of monomials up to a finite degree. The conjugacy $\hv$ could then be used to identify a nonlinear change of variable, under which the iterates of the latent variables $\yv_n := \hv(\xv_n)$ are governed by a polynomial map. A well-known example of this is the conjugacy between the tent map
\begin{equation}
	f(x) = 2\min\{x,1-x\}
\end{equation}  
and the logistic map
\begin{equation}
	g(y) = 4y(1-y).
\end{equation}
Indeed, the conjugacy is known to be $h(x) = \frac{2}{\pi}\arcsin(\sqrt{x})$, thus providing a nonlinear change of coordinates that takes iterates of a piecewise continuous mapping, $f$, to a quadratic polynomial map, $g$.  

As mentioned above, an advantage to using the neural network to obtain conjugate mappings is that it provides a notion of dimensionality reduction. Recall that the conjugacy is only between attractors and the dimension of the attractor is approximated by the Kaplan--Yorke dimension, $D_{KY}$. If the flow of the governing system lies in $\mathbb{R}^D$, then necessarily we have that $D_{KY} \leq D$. In almost every case of interest this inequality is strict. Recalling from Subsection~\ref{subsec:LEs} that we take the latent mapping $\gv$ to have phase-space $\mathbb{R}^d$ with $d \approx D_{KY} - 1$, we have the potential to greatly reduce the dimension from that of the observations $\{\xv_n\} \subset \mathbb{R}^{D-1}$ to that of the latent variables $\yv \in \mathbb{R}^d$. This dimensionality reduction provides one of the major benefits for obtaining conjugacies with a neural network since it allows us to understand relatively low-dimensional attractors in high-dimensional and even infinite-dimensional systems. We will illustrate this utility with the Kuramoto--Sivashinsky equation in Section~\ref{sec:KS}.  

The loss function used to train the network is the sum of the different losses:
\begin{enumerate}
	\item {\bf Loss 1: Conjugacy loss.} A requirement for obtaining a conjugacy is that $\hv$ be a homeomorphism, and is therefore invertible. To ensure this, we require $\hv^{-1} \circ \hv$ to reconstruct the training data. This loss is given by $\|\xv_n - \hv^{-1}(\hv(\xv_n))\|_\mathrm{MSE}$, where $\|\cdot\|_\mathrm{MSE}$ refers to the mean square error.
	\item {\bf Loss 2: Prediction loss.} Applying the network to $\xv_n$ results in $\hv^{-1}(\gv(\hv(\xv_n)))$, and to ensure this replicates the dynamics of the unknown Poincar\'e map $\fv$ we want the output to be $\xv_{n+1} = \fv(\xv_n)$. To enforce this we use the loss $\|\xv_{n+1} - \hv^{-1}(\gv(\hv(\xv_n)))\|_\mathrm{MSE}$. Furthermore, a conjugacy maps orbits of $\fv$ to orbits of $\gv$, and so we may generalize the loss to incorporate this. In general we have 
	\begin{equation}
		\frac{1}{s}\sum_{j = 1}^s \|\xv_{n+j} - \hv^{-1}(\gv^j(\hv(\xv_n)))\|_\mathrm{MSE} 
	\end{equation}
	where the number of steps $s \geq 1$ is a hyperparameter in the system and $\gv^j$ is the composition of $\gv$ with itself $j$ times.
	\item {\bf Loss 3: Latent mapping loss.} Starting with $\xv_n$ and $\xv_{n+1}$, following the commutative diagram \eqref{ComDiag} means that $\gv(\hv(\xv_n))$ and $\hv^{-1}(\xv_{n+1})$ should be the same. This is reflected by the loss $\|\hv^{-1}(\xv_{n+1}) - \gv(\hv(\xv_n))\|_\mathrm{MSE}$, which can similarly be generalized to $s\geq 1$ steps into the future with
	\begin{equation}
		\frac{1}{s}\sum_{j = 1}^s \|\hv^{-1}(\xv_{n+j}) - \gv^j(\hv(\xv_n))\|_\mathrm{MSE}
	\end{equation}  
	\item {\bf Loss 4: Network regularization.} \new{We apply elastic net regularization by including $\ell^1$ and $\ell_2$ regularization terms to avoid overfitting.}  
\end{enumerate} 
In some cases we have found it advantageous to introduce one more loss:
\begin{enumerate}
	\item[5.] {\bf Loss 5: Attractor stretching.} If the image of $\hv$ is very small then Loss 3 will also be small, even if the latent mapping $\gv$ is inaccurate. To get around this we introduce a loss that stretches the attractor across the hypercube $[0,1]^d$. This loss takes the form
	\begin{equation}
		\frac{1}{4d}\sum_{j = 1}^d (1 - \max_n\{h_j(\xv_n)\})^2 + (\min_n\{h_j(\xv_n)\})^2 + (1 - \max_n\{g_j(\hv(\xv_n))\})^2 + (\min_n\{g_j(\hv(\xv_n))\})^2, 
	\end{equation} 
	where $h_j$ and $g_j$ denotes the $j$th component of $\hv$ and $\gv$, respectively. This loss guarantees that the domain and range of $\gv$ goes from $0$ to $1$ in each component of $\hv(\xv_n)$. Note that such a stretching or compressing in each dimension is achieved by conjugating any map using a linear conjugacy function. Therefore, this loss works to make the choice of $\gv$ unique by fixing the the extremal values of the domain and range. 
\end{enumerate}  
\new{Network weights are initialized randomly according to a normal distribution with mean zero, and each hidden layer of the network uses the SeLu activation function~\cite{Selu}. To optimize performance, all training data is scaled to only take values between 0 and 1. Complete python code and Jupyter notebooks to implement these methods are available at {\bf github.com/jbramburger/Deep-Conjugacies}, along with saved models that reproduce the results in this work. These saved models are the result of extensive searches for the smallest loss over the various parameters used to build the network and the above components of the loss function. These parameter values are given in the appendix.}

\section{Low-Dimensional Systems}\label{sec:Low} 

We begin our investigation with three chaotic systems with three-dimensional phase spaces. These examples are meant to illustrate the various tasks that the network is able to perform. Particularly, we will see that in the case of the R\"ossler system the network can perform a nonlinear change of coordinates that improves the fit of the data to a polynomial Poincar\'e map. In the case of the Lorenz system we will see how the network can be used to obtain conjugacies between the section data and heuristic mappings that have long been used to understand the section dynamics. Finally, we will use the Gissinger system to demonstrate the dimensionality reduction aspect of the network by moving from the three-dimensional phase space of the continuous flow, to the two-dimensional Poincar\'e section data, to a conjugated one-dimensional mapping.

\subsection{R\"ossler System} 

\begin{figure} 
\center
\includegraphics[width = 0.35\textwidth]{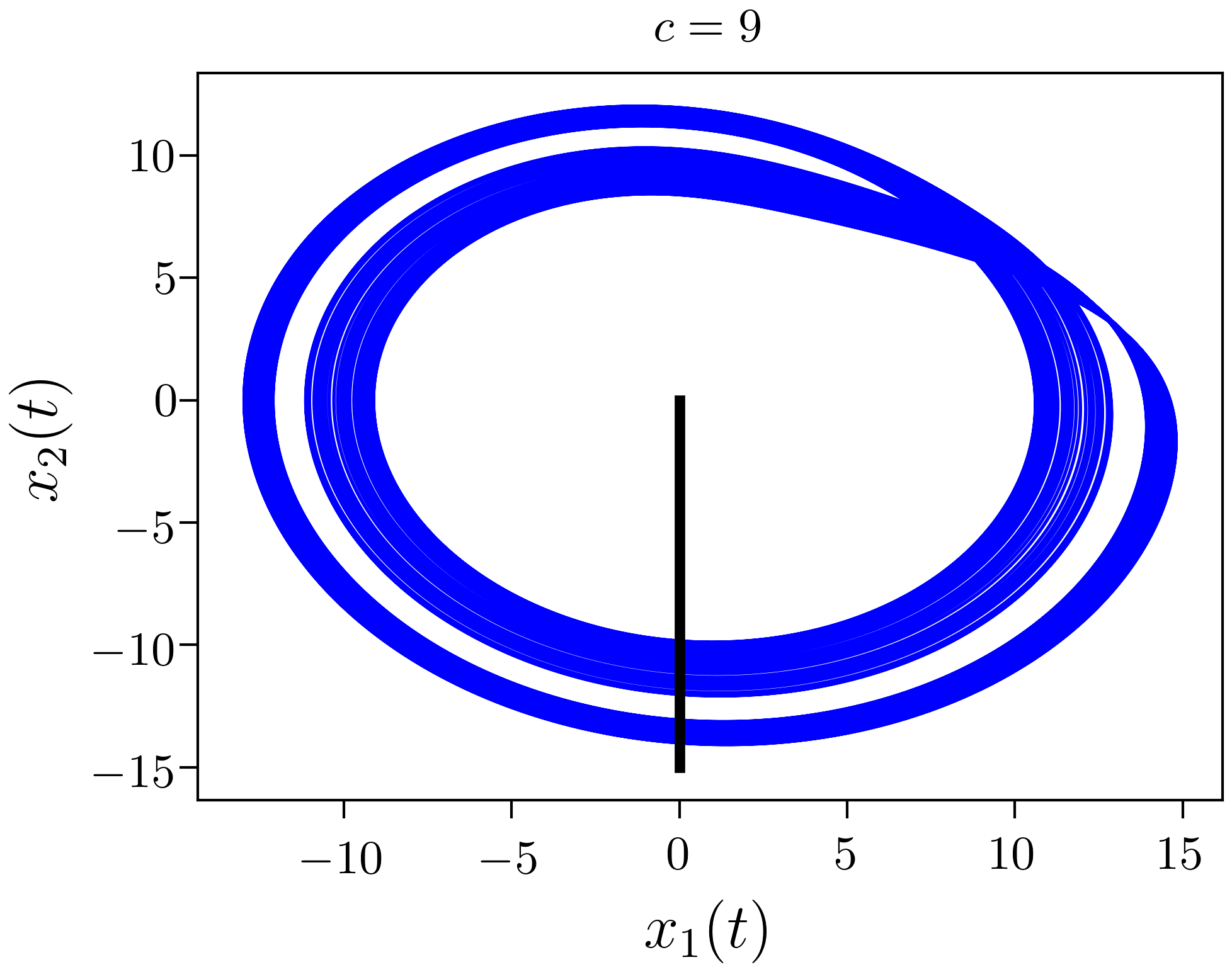}  \quad \includegraphics[width = 0.35\textwidth]{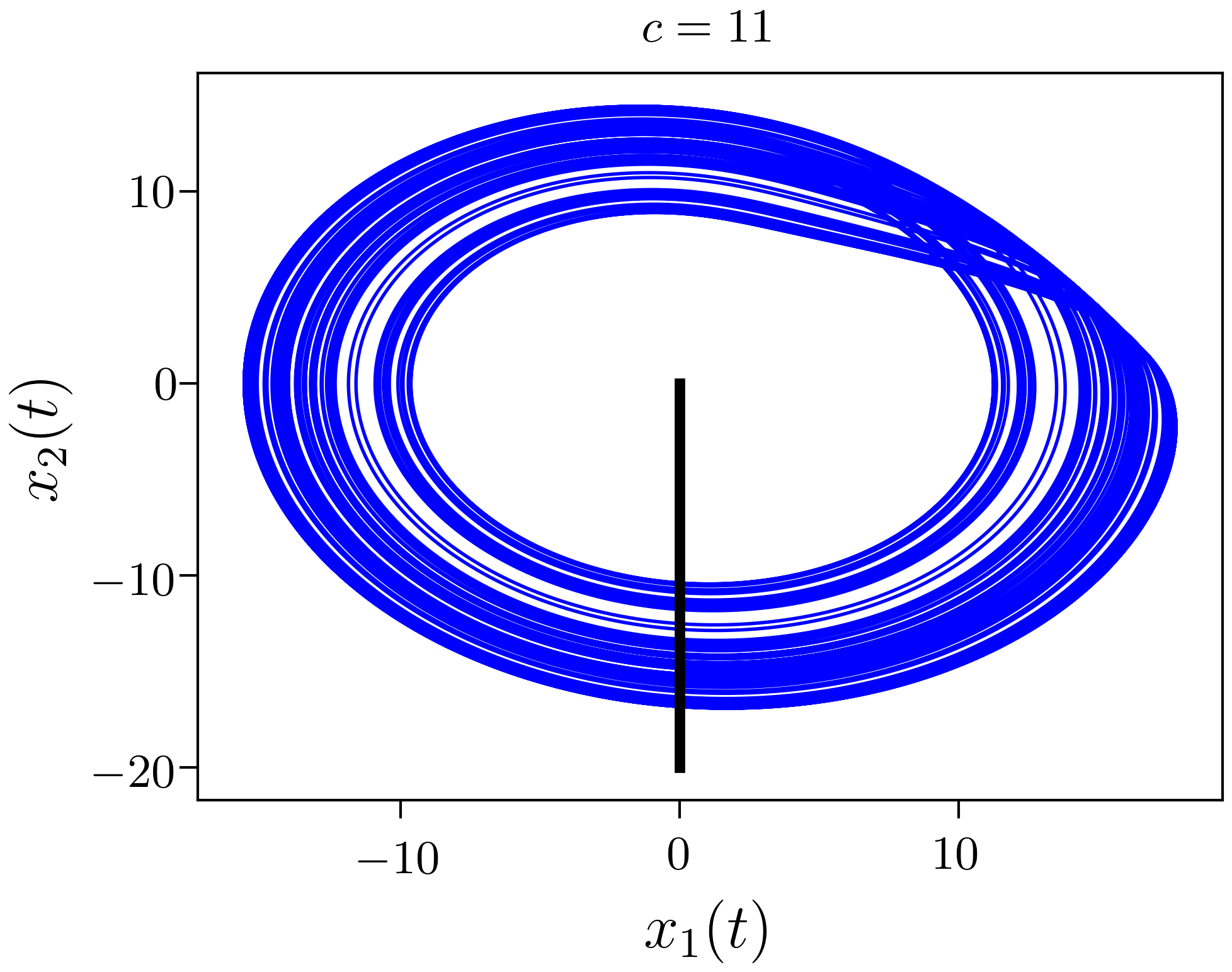} \\
\includegraphics[width = 0.35\textwidth]{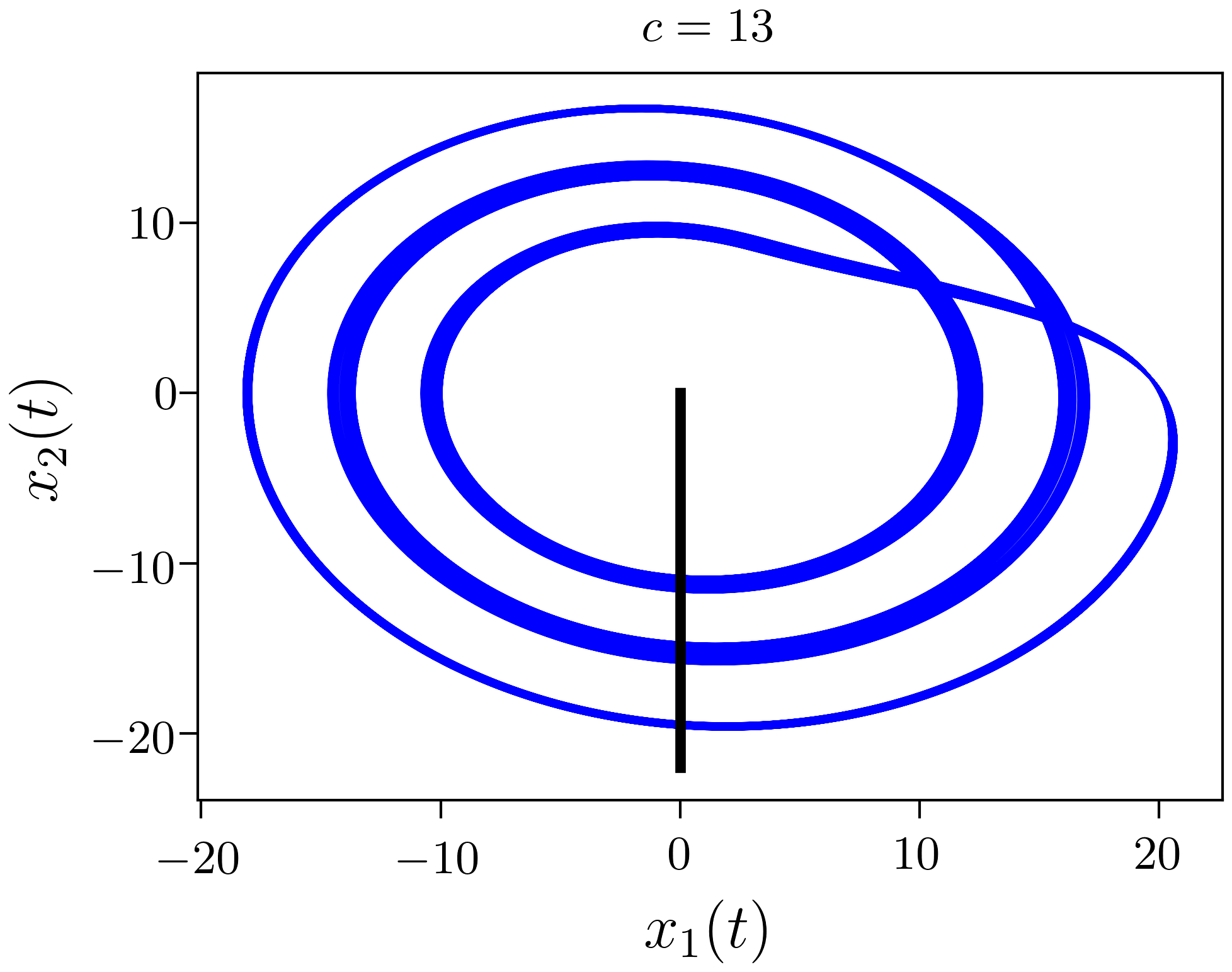}  \quad \includegraphics[width = 0.35\textwidth]{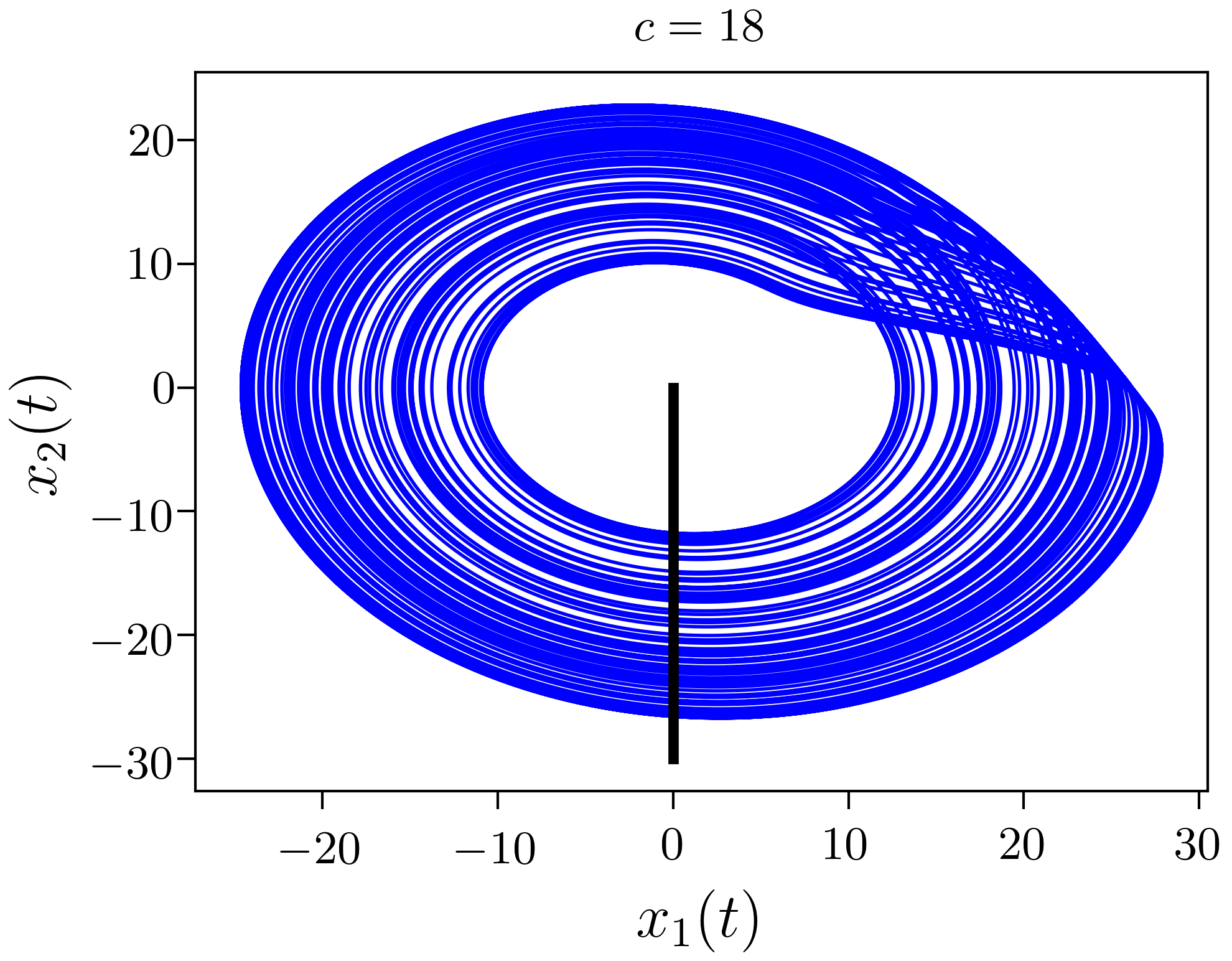} 
\caption{Typical trajectories of \eqref{Rossler} on the chaotic attractor for varying values of $c$ in the $(x_1,x_2)$-plane. The black line is a projection of the Poincar\'e section $x_1 = 0$ and $\dot{x}_1 > 0$.}
\label{fig:Rossler}
\end{figure} 

\begin{figure} 
\center
\includegraphics[width = 0.35\textwidth]{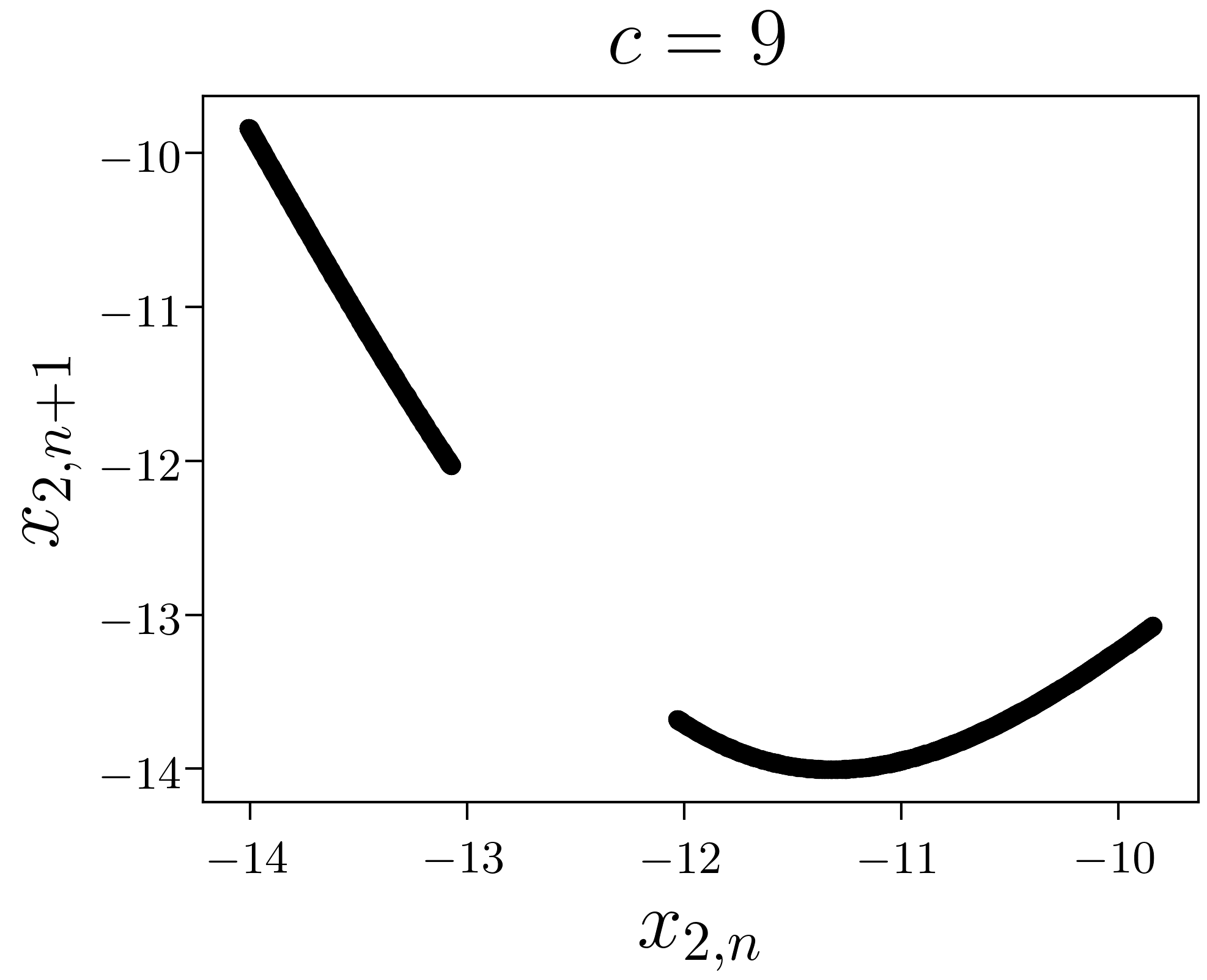}  \quad \includegraphics[width = 0.35\textwidth]{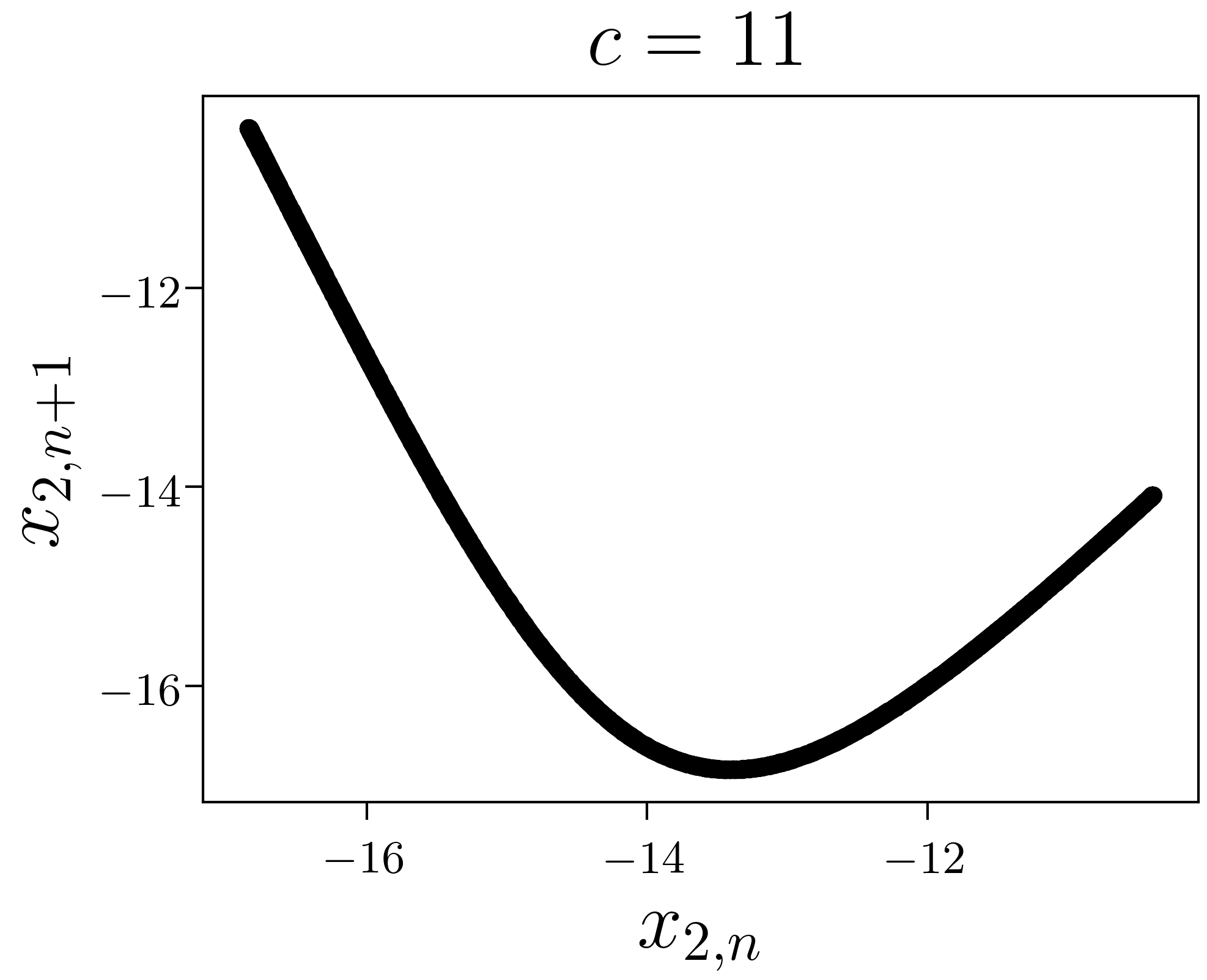} \\
\includegraphics[width = 0.35\textwidth]{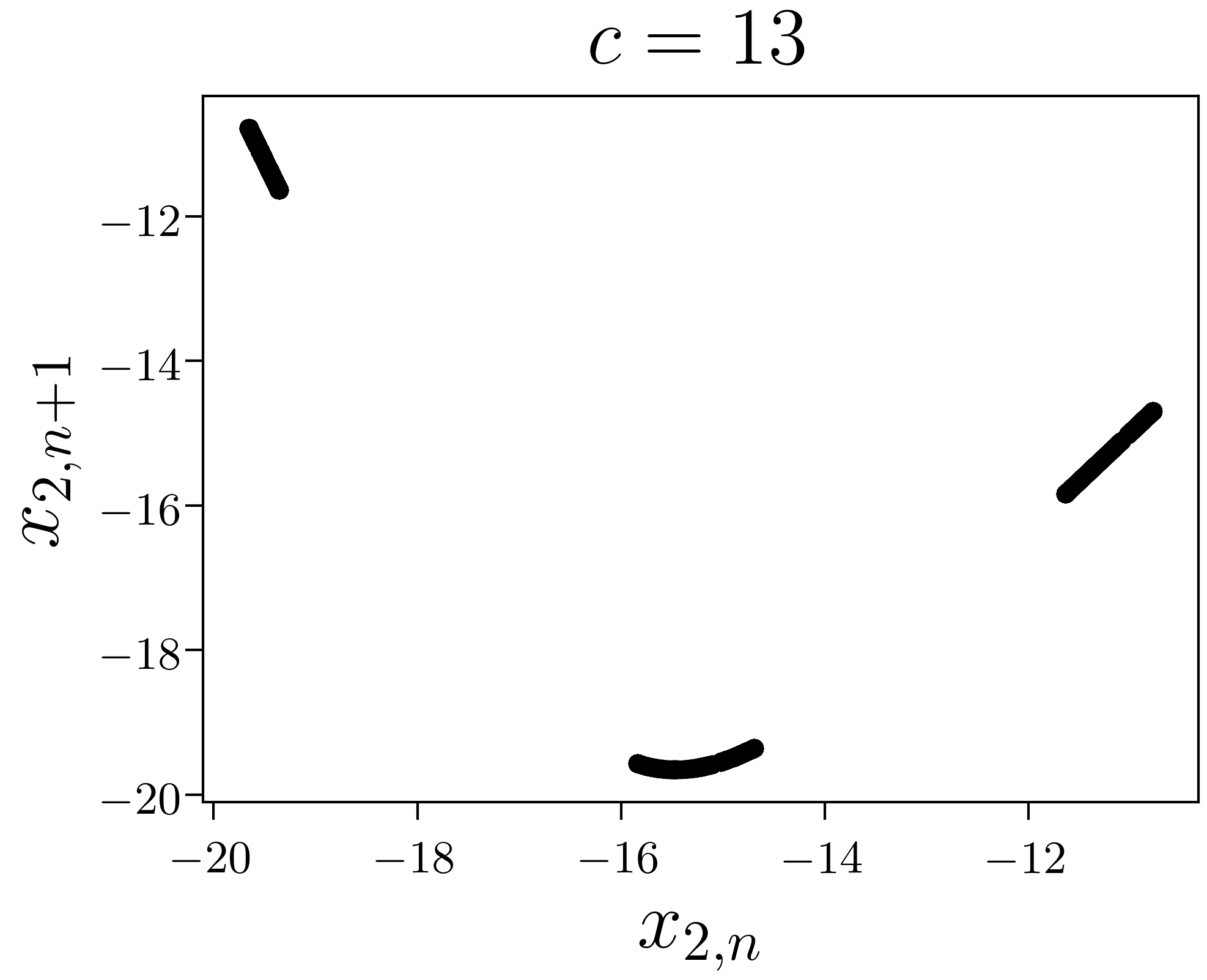}  \quad \includegraphics[width = 0.35\textwidth]{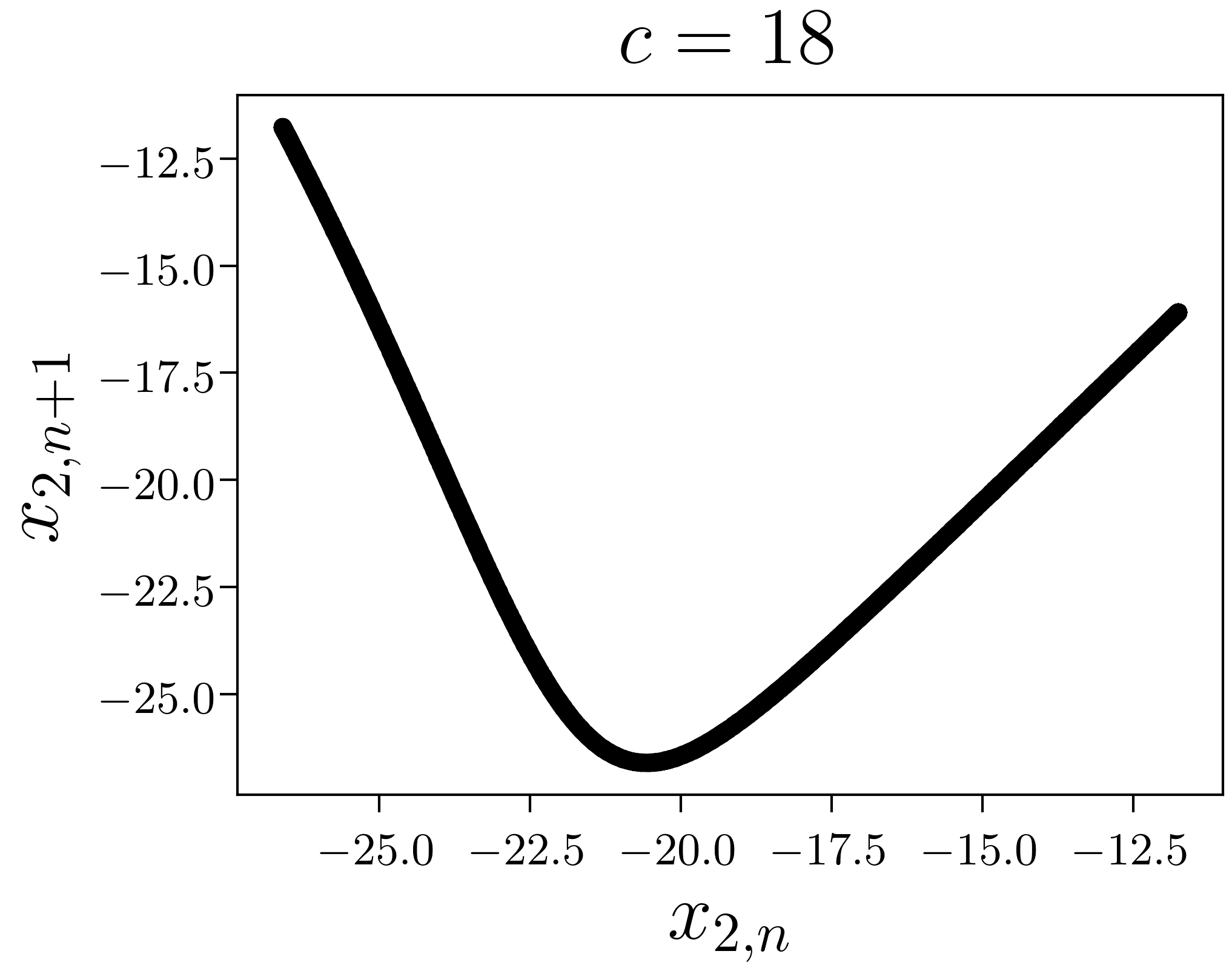} 
\caption{First return maps on the attractor of the R\"ossler system \eqref{Rossler} for varying $c$.}
\label{fig:RosslerSection}
\end{figure}

Our first case study will be the R\"ossler system, given by
\begin{equation}\label{Rossler} 
	\begin{split}
		\dot{x}_1 &= -x_2 - x_3 \\
		\dot{x}_2 &= x_1 + 0.1x_2 \\
		\dot{x}_3 &= 0.1 + x_3(x_1-c)
	\end{split}
\end{equation}
where $c\in\mathbb{R}$ is a bifurcation parameter. It is well known when $x_1$ crosses $0$ on the attractor, so does $x_3$. Thus, we take $x_1 = 0$ with $\dot x_1 > 0$ as our Poincar\'e section, reducing ourselves to iterating in a single variable $x_2$ as we cross the section. In Figure~\ref{fig:RosslerSection} we provide plots of the first return map on the attractor for varying $c$. We observe that the iterates appear to lie approximately along a quadratic function, prompting many to use the logistic map 
\begin{equation}\label{Logistic}
	y_{n+1} = ry_n(1-y_n), \quad r \in (0,4]
\end{equation}
as a simple heuristic for understanding the dynamics of \eqref{Rossler} on the attractor.  

Previous attempts to fit the first returns maps to simple polynomial functions were met with varying degrees of success~\cite{Bramburger}. For parameter values $c$ where the attractor intersects the Poincar\'e section at finitely many points (non-chaotic), the SINDy method~\cite{SINDy} is able to obtain polynomial mappings that accurately describe these iterations. This is not surprising since we are only required to fit the polynomial to finitely many points. As the $c$ value increases and the attractor becomes chaotic, the polynomial functions become less well-suited to the first return maps of the R\"ossler system \eqref{Rossler}. The inability for polynomial mappings to faithfully describe the chaotic return map dynamics thus limits our ability to forecast the chaotic system, understand the statistics of the attractor, and pull out periodic orbits~\cite{BramburgerUPO}.  

One method to obtain a function that better represents the first return map of the R\"ossler attractor is to include more than just polynomials in the library of candidate functions for which the first return map must belong to the span of. There is an immediate limitation to this though, and that is deciding which functions to include. Instead, we will use the network to obtain a nonlinear change of variables which conjugates the dynamics of the first return map of the R\"ossler attractor into a simple quadratic function of the form 
\begin{equation}\label{RossQuad}
	g(y) = c_1y + c_2y^2,
\end{equation}
where $c_1,c_2$ are variables to be learned by the network. We first note that we do not include any constant term in our mapping $g$ since a simple linear shift of the variable $y$ can eliminate this. Therefore, we leave this linear shift to be learned by the network. Second, a linear stretching/compressing of the variable $y$ will conjugate the map $\eqref{RossQuad}$ into the map
\begin{equation}
	\tilde{g}(y) = c_1y\bigg(1 + \frac{c_2}{|c_2|}y\bigg),
\end{equation} 
and therefore if $c_2 < 0$ we get that the mapping \eqref{RossQuad} is dynamically equivalent to the logistic map with $r = c_1$. The election to use two variables, $c_1$ and $c_2$, allows the network greater freedom to fine-tune the conjugacy between the training data and the iterates of the latent map.  

\begin{table}
\centering
\renewcommand{\arraystretch}{1.2}
\begin{tabular}{|c|c|c|c|c|}
\hline
{}  &  \multicolumn{2}{c|}{Neural Network} & \multicolumn{2}{c|}{SINDy Method}\\
\hline
Parameter  & Mapping   & $r$-value    & Mapping   & $r$-value\\
\hline
c = 9   &  $3.6075y -4.9044y^2 $ & 3.6075   & $1.0219 - 2.9885y + 2.2063y^2$  & 3.6219\\
c = 11   &  $3.8302 y -4.4801y^2$ & 3.8302   & $1.0217 - 3.4481y + 2.9242y^2$  & 3.7985\\
c = 13  &  $3.8502y -6.5346y^2$  &  3.8502   & $1.0075 - 3.5638y + 3.1613y^2$  & 3.8456 \\
c = 18  &  $3.9661y - 4.6718y^2$  &  3.9661   & $0.9776 -3.5219 y + 3.4089y^2$  & 3.6737\\
\hline
\end{tabular}
\caption{Comparison of discovered R\"ossler Poincar\'e mappings along with the corresponding value of $r$ in the Logistic map~\eqref{Logistic} that they are conjugate to.}
\label{tbl:Rossler}
\end{table}

\begin{figure} 
\center
\includegraphics[width = 0.45\textwidth]{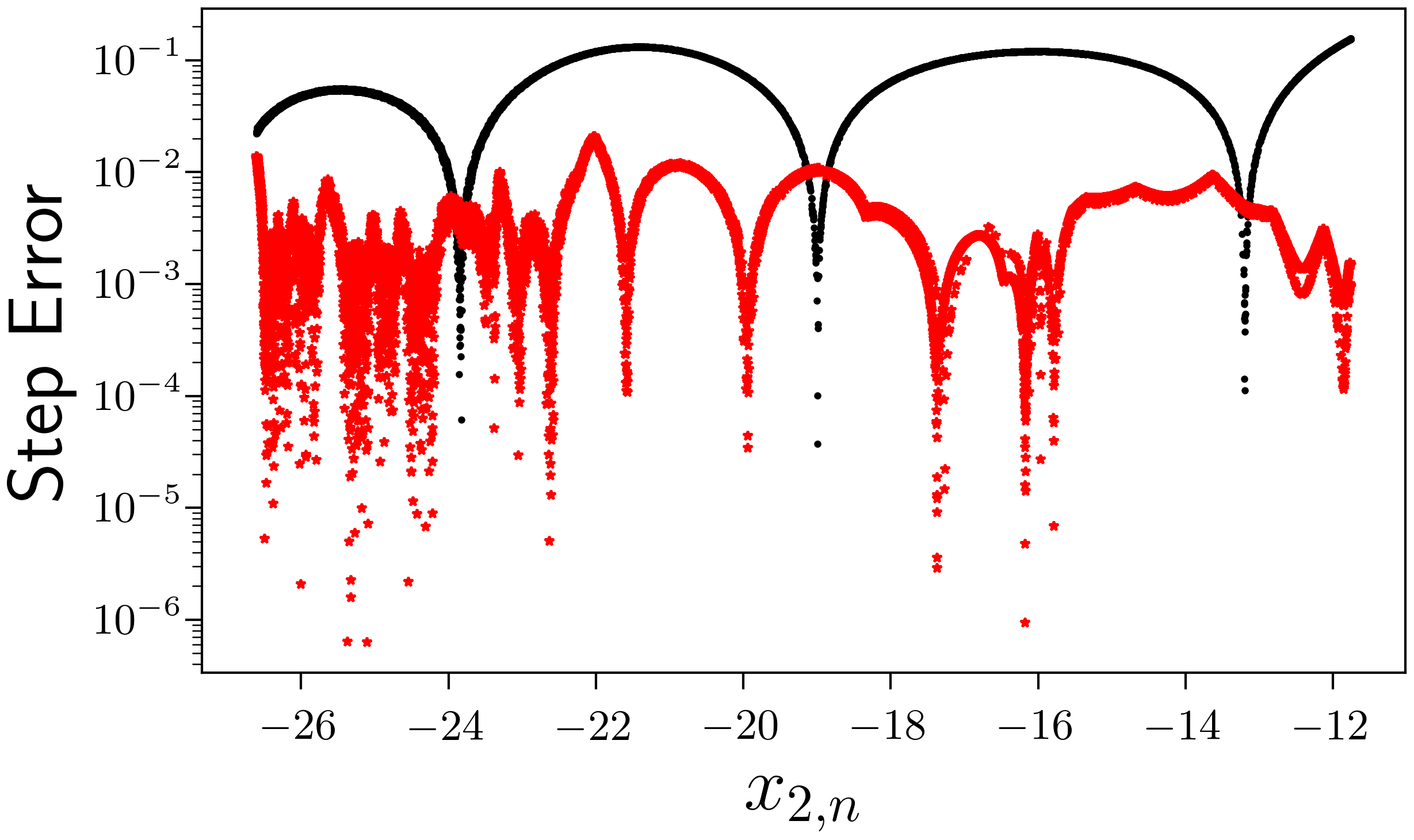} \ \includegraphics[width = 0.45\textwidth]{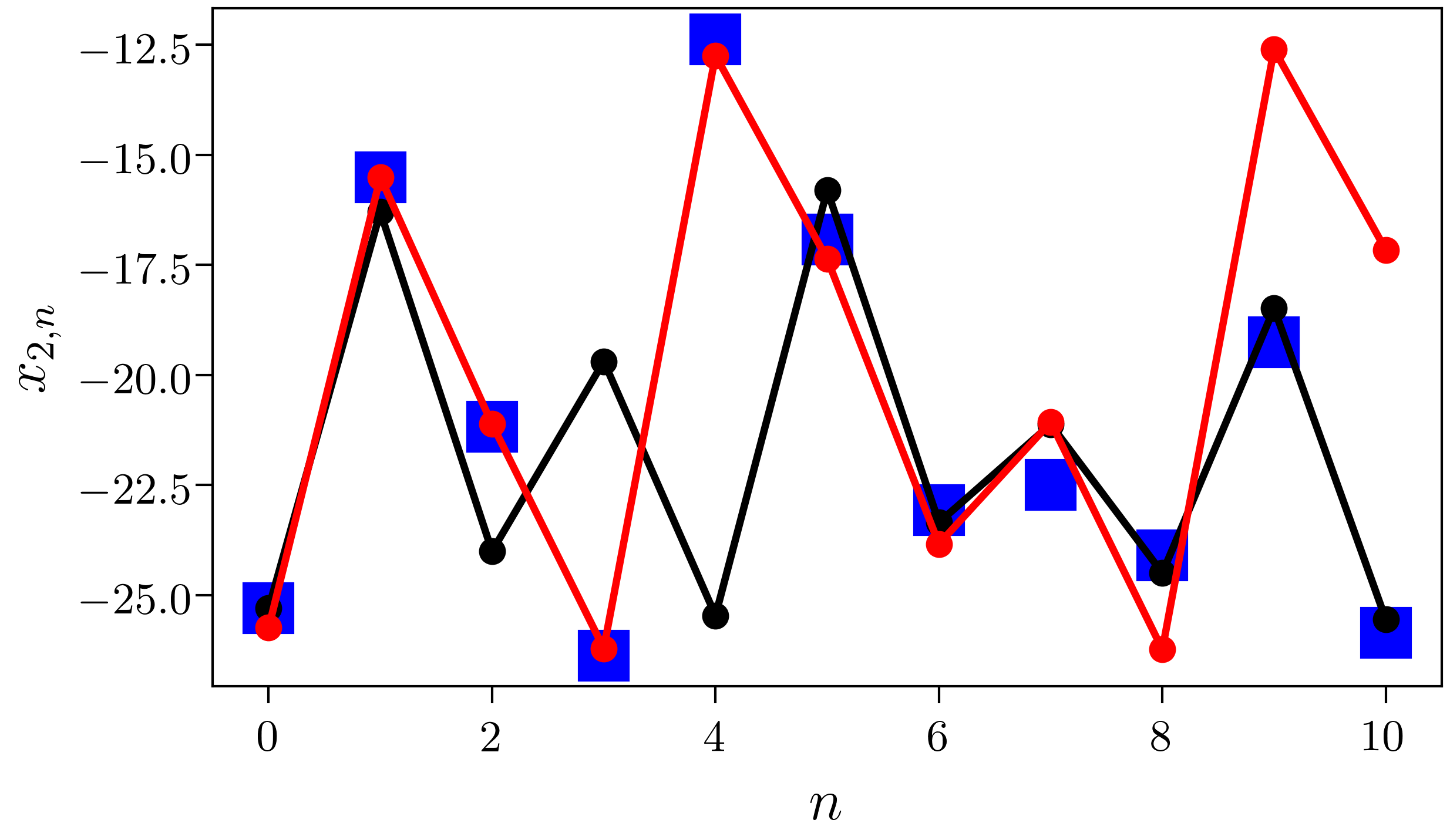} 
\caption{Left: One-step error of every training data point on the attractor of the R\"ossler system with $c = 18$. Plotted is a comparison between the discovered SINDy model (black, circles) and neural network (red, stars). Right: The reduced step error with the neural network results in more accurate forecasting in the Poincar\'e section. Compare the training data (blue, squares) with the SINDy model forecast (black) and the neural network model forecast (red), all with the same random initial condition on the attractor.}
\label{fig:RosslerError}
\end{figure}

The results for $c = 9,11,13$ and $18$ are presented in table~\ref{tbl:Rossler}. We also provide the quadratic function produced by the SINDy method using the same training data. The obtained quadratic SINDy mappings are easily conjugated to the logistic map \eqref{Logistic}, and so we also provide the corresponding logistic parameter $r$ for comparison. We have the best agreement between the discovered models at $c = 9,13$, when the attractor and the LEs are relatively small. The greatest discrepancy is at $c = 18$, where the attractor and the LEs are relatively large. To compare the results at $c = 18$, in Figure~\ref{fig:RosslerError} we provide the one-step error between the two maps and the training data. That is, we iterate every training data point on the attractor through the network (red) and through the SINDy model (black) and measure the distance from the true value, as obtained through numerically integrating the R\"ossler system. This relatively large error in the SINDy model means that attempts to forecast the iterates in the section will become inaccurate quickly due to the sensitivity of the mapping to both initial conditions and the coefficients. We illustrate this inaccuracy in the right image of Figure~\ref{fig:RosslerError} with a random initial condition iterated forward using the network and the SINDy method. We also plot the true iterates as obtained through numerical integration as blue squares for comparison. We see that the SINDy model becomes inaccurate after about two steps, while the network mapping doesn't become inaccurate until about seven steps. Finally, we emphasize that between iterates in the section is approximately 7 time units, meaning the network has forecast the trajectory nearly 50 time units into the future.     

We further highlight the utility of the conjugate mapping in the context of controlling chaos. In a related work it was shown that the SINDy models in chaotic regions are able to locate period one and two UPOs~\cite{BramburgerUPO}, but attempts to located higher period UPOs was prevented due to numerical inaccuracy. To demonstrate the utility of our results, we report here that using the latent mapping $g$ we are able to locate and stabilize periodic orbits at $c = 11$ up to at least period six. These results are not presented here for brevity, but these calculations can be reproduced using an accompanying MATLAB script included in the code repository for this article. We do highlight the fact that these orbits were stabilized using the methods of Subsection~\ref{subsec:Control} and the derivatives were estimated directly from numerical integrations of the full R\"ossler system. 

\new{Finally, we remind the reader that the trained neural networks are the result of extensive hyperparameter searches. The values for the hyperparameters for each value of $c$ is reported in the appendix. At present we do not have any analytical results on the robustness of the discovered mappings with respect to the neural network hyperparameters, but preliminary numerical explorations have been promising. For example, taking $c = 11$ and varying the hyperparameters in a neighbourhood of the values reported in the appendix returns similar values of the loss function. Particularly, taking sufficiently widths $\geq 100$ and one or more blocks in and out in the network results in nearly identical results to those reported here. Therefore, we hope to provide convergence guarantees and robustness results in a follow-up investigation.}

\subsection{Lorenz System} 

\begin{figure} 
\center
\includegraphics[width = 0.35\textwidth]{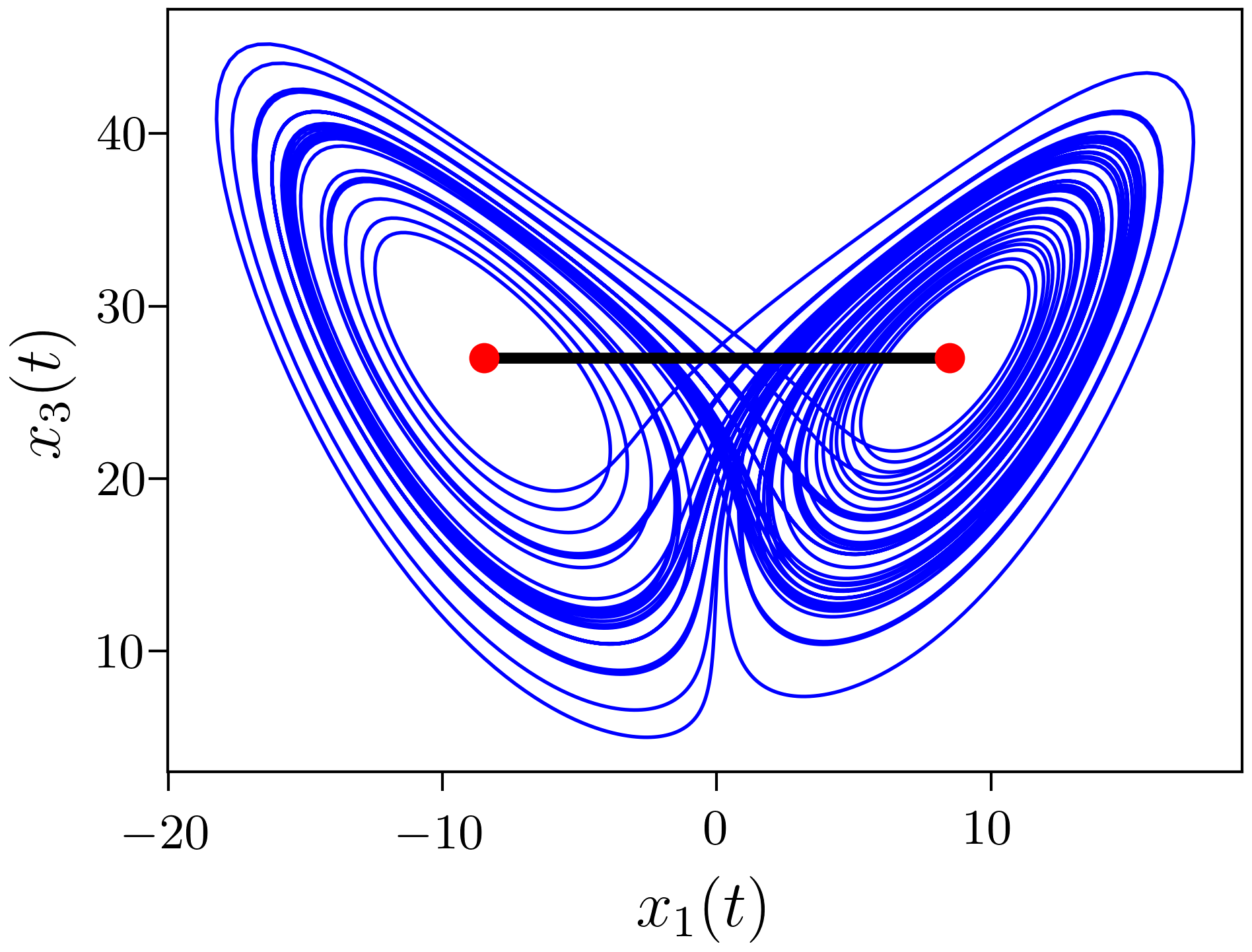} \ \includegraphics[width = 0.35\textwidth]{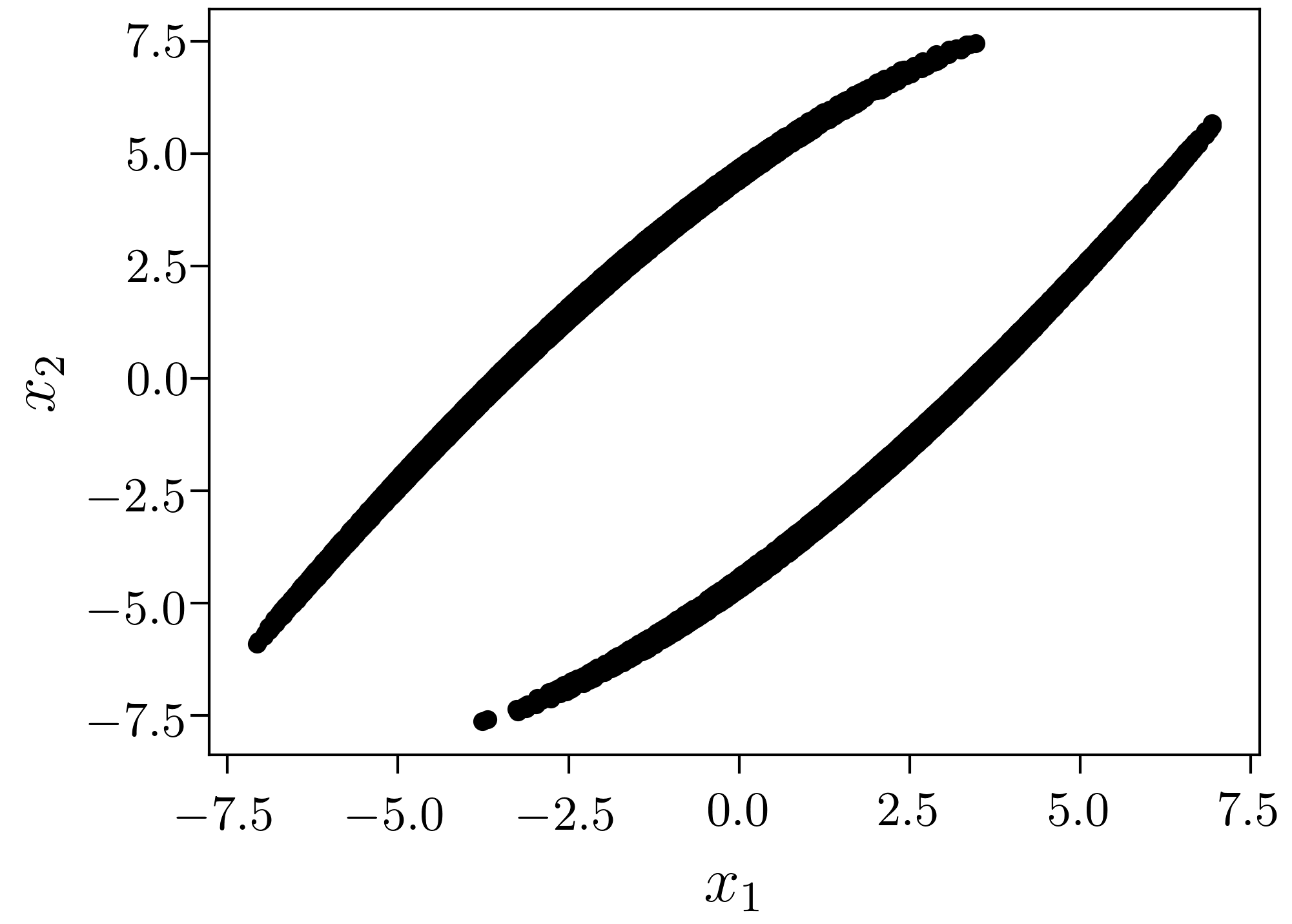}
\caption{Left: A typical chaotic trajectory of the Lorenz system \eqref{Lorenz} projected into the $(x_1,x_3)$-plane. In red are the two nontrivial equilibria and black line represents a side view Poincar\'e section given by $x_3 = 27$ and $\dot x_3 < 0$. Right: The intersection of the chaotic attractor with the Poincar\'e section.}
\label{fig:Lorenz}
\end{figure} 

Let us turn our attention now to one of the most celebrated equations in all of chaotic dynamics: the Lorenz system. The governing equations are given by
\begin{equation}\label{Lorenz} 
	\begin{split}
		\dot{x}_1 &= 10(x_2 - x_1) \\
		\dot{x}_2 &= x_1(28 - x_3) - x_2 \\
		\dot{x}_3 &= x_1x_2 - \frac{8}{3} x_3
	\end{split}
\end{equation}
where we are using the standard parameter values that induce chaotic dynamics. We take the Poincar\'e section to be $x_3 = 27$ and $\dot{x}_3 < 0$, which contains the two nontrivial equilibria at the centres of the chaotic lobes. In Figure~\ref{fig:Lorenz} we plot a typical chaotic trajectory of system \eqref{Lorenz}, the two nontrivial equilibria, and the Poincar\'e section in the $(x_1,x_3)$-plane. The classical method to understand the Poincar\'e map is to follow~\cite{Guckenheimer}, where they assume there is a change of variable from the $(x_1,x_2)$ coordinate system in the section to a new set of coordinates $(y_1,y_2)$ for which the Poincar\'e map can be written in the skew-product form
\begin{equation}\label{LorenzMap}
	\gv(y_1,y_2) = \begin{bmatrix}g_1(y_1) \\ g(y_1,y_2)\end{bmatrix}.
\end{equation}   
Furthermore, the map is assumed to satisfy $\gv(-y_1,-y_2) = -\gv(y_1,y_2)$ and has a jump discontinuity at $y_1 = 0$.

Despite having no exact correspondence between the Poincar\'e map of the Lorenz system and the function \eqref{LorenzMap}, it has primarily been through such skew-product maps that we have formed our modern understanding of the chaotic Lorenz system. To strengthen this connection, we attempt to learn a conjugacy between the training data shown on the right of Figure~\ref{fig:Lorenz} and a mapping of the form \eqref{LorenzMap} using our proposed neural network architecture. We have elected to attempt to conjugate this data with a skew-product map of the form
\begin{equation}\label{LorenzMapNN}
	\begin{split}
		g_1(y_1) &= -\mathrm{sgn}(y_1) + c_1y_1 + c_2|y_1|y_1 \\
		g_2(y_1,y_2) &= d_0\mathrm{sgn}(y_1) + d_1y_2,
	\end{split}
\end{equation}
where $\mathrm{sgn}(\cdot)$ is the function that returns 1 if the argument is positive and $-1$ if it is negative. The specific form of $g_2$ was proposed in~\cite{LorenzSymbols}, thus motivating the choice here. The form for $g_1$ is taken to be the simplest quadratic model that satisfies the constraints assumed in~\cite{Guckenheimer}. The goal of employing the neural network is not only to find the coordinate transformation from $(x_1,x_2)$ to $(y_1,y_2)$, but to obtain an invertible nonlinear change of coordinates that simplifies the right-hand-side of the mapping as well.

\begin{table} 
\centering
\renewcommand{\arraystretch}{1}
\begin{tabular}{ |c||c|c|c| }
\hline
Symbolic Sequence & Period & fsolve Iterates & fsolve Time (seconds)  \\
\hline
LR & 1.5560 & 7 & 0.9253 \\
LLR & 2.3032 & 8 & 2.8569 \\
LLLR & 3.0208 & 9 & 7.8822  \\
LLRR & 3.0816 & 7 & 7.3054 \\
LLLLR & 3.7228 & 8 & 12.2913 \\
LLLRR & 3.8174 & 9 & 20.1933 \\
LLRLR & 3.8667 & 8 & 12.0158 \\
\hline
\end{tabular}
\caption{Symbolic sequences of periodic orbits in the discovered conjugate mapping are used to obtain periodic orbits in the Lorenz system \eqref{Lorenz}. All symbolic sequences up to symmetry up to length 5 are presented, along with the period of the related UPO to the Lorenz system, the number of iterates to converge to the UPO, and the time taken.}
\label{tbl:Lorenz}
\end{table} 

Using $26631$ training data points gathered from a single trajectory of system \eqref{Lorenz}, we find a numerical conjugacy to a mapping of the form \eqref{LorenzMapNN} with 
\begin{equation}\label{LorenzCoefficients}
	\begin{split}
		c_1 &= 2.5248 \quad \quad d_0 = -0.34275 \\
		c_2 &= 1.6595 \quad \quad d_1 = 1.7825.
	\end{split}
\end{equation} 
\new{The conjugacy to the mapping \eqref{LorenzMapNN} with the coefficients \eqref{LorenzCoefficients} now allows one to extract the intersection of the UPOs of \eqref{Lorenz} with the Poincar\'e section. From the method outlined in Section~\ref{subsec:UPO}, this in turn allows one to seed initial guesses for the UPOs to be found in the continuous time flow. For example, a period 2 point of the map \eqref{LorenzMapNN} is obtained by solving 
\begin{equation}\label{LorenzPer2}
	\begin{split}
		y_1 &= g_1(g_1(y_1)) \\
		y_2 &= g_2(g_1(y_1),g_2(y_1,y_2))
	\end{split}
\end{equation}
for $(y_1,y_2)$. Indeed, such a point is mapped back to itself after exactly two iterations of \eqref{LorenzMapNN}. Denoting the extracted period 2 point that solves \eqref{LorenzPer2} by $(\bar y_1,\bar y_2)$, we then use the discovered function $\hv^{-1}$ from the neural network to define $\bar \xv_1 = \hv^{-1}(\bar y_1, \bar y_2)$ and $\bar \xv_2 = \hv^{-1}(g_1(\bar y_1), g_2(\bar y_1,\bar y_2))$, using the notation of Section~\ref{subsec:UPO}. Since periodic points are mapped into each other by the homeomorphism $\hv$, it follows that $\bar\xv_1$ and $\bar\xv_2$ are period 2 points in the Poincar\'e section of \eqref{Lorenz}, thus belonging to a UPO that crosses the section exactly twice before completing a full period. With these points we can then seed initial guesses for the continuous-time UPOs of \eqref{Lorenz} using the method outlined in Section~\ref{subsec:UPO}.}

\new{To demonstrate the utility of the discovered mapping, we will follow the above outlined procedure to obtain UPOs of \eqref{Lorenz}. We use the built-in MATALB function fsolve to obtain periodic points of the map \eqref{LorenzMapNN} and to solve \eqref{UPOFinder} for the desired continuous-time UPOs. All tolerances have been set to at least $10^{-15}$.} Our results are summarized in table~\ref{tbl:Lorenz} where we have enumerate information for all periodic orbits that have up to 5 intersections with the Poincar\'e section before getting back to where they started. The `L' denotes the trajectory winding around the left lobe of the attractor, while the `R' denotes the trajectory winding around the right lobe. Note that system \eqref{Lorenz} is equivariant with respect to the action $(x_1,x_2,x_3) \mapsto (-x_1,-x_2,x_3)$, meaning that a sequence of L's and R's representing a periodic orbit is mapped by the symmetry to another periodic orbit obtained by flipping the L's to R's and the R's to L's. Hence the sequence LRR, for example, is absent from table~\ref{tbl:Lorenz} since the corresponding UPO can be obtained from the sequence LLR and therefore has the same period. The periods of the UPOs agree with those in~\cite{LorenzSymbols} up to at least three significant digits, while greater accuracy can be obtained by increasing the temporal resolution. We note that increasing the resolution does increase the time taken for each iterate, but rarely does it increase the number of iterates. This testifies to the accuracy of the conjugacy between the Poincar\'e map of the Lorenz system and to our discovered mapping. The corresponding UPOs can be plotted using a simple MATLAB script, included in the code repository for this article.  

The skew-product form of the discovered mapping \eqref{LorenzMapNN} emphasizes that almost all of the information about the attractor in the section can be understood through the dynamics of the $y_1$ variable. This is further emphasized by examining the Lyapunov spectrum. The LEs are given by $\lambda_1 = 0.90$, $\lambda_2 = 0$, and $\lambda_3 = -14.57$, giving a Kaplan--Yorke dimension of $2.06$, and implying that a one-dimensional map should be able to accurately describe the dynamics. The projection of the attractor of our discovered mapping $\gv(y_1,y_2)$ onto the first component gives a bijection since the periodic orbits of $g_1(y_1)$ lie in one-to-one correspondence with the periodic orbits of $\gv(y_1,y_2)$ which densely fill the attractor. With these facts as motivation, we are thus able to use the network to obtain a conjugacy between the Lorenz section data and the mapping $g_1(y_1)$, providing for a successful dimensionality reduction of the dynamics on the attractor.

\subsection{Gissinger's System} 

\begin{figure} 
\center
\includegraphics[width = 0.35\textwidth]{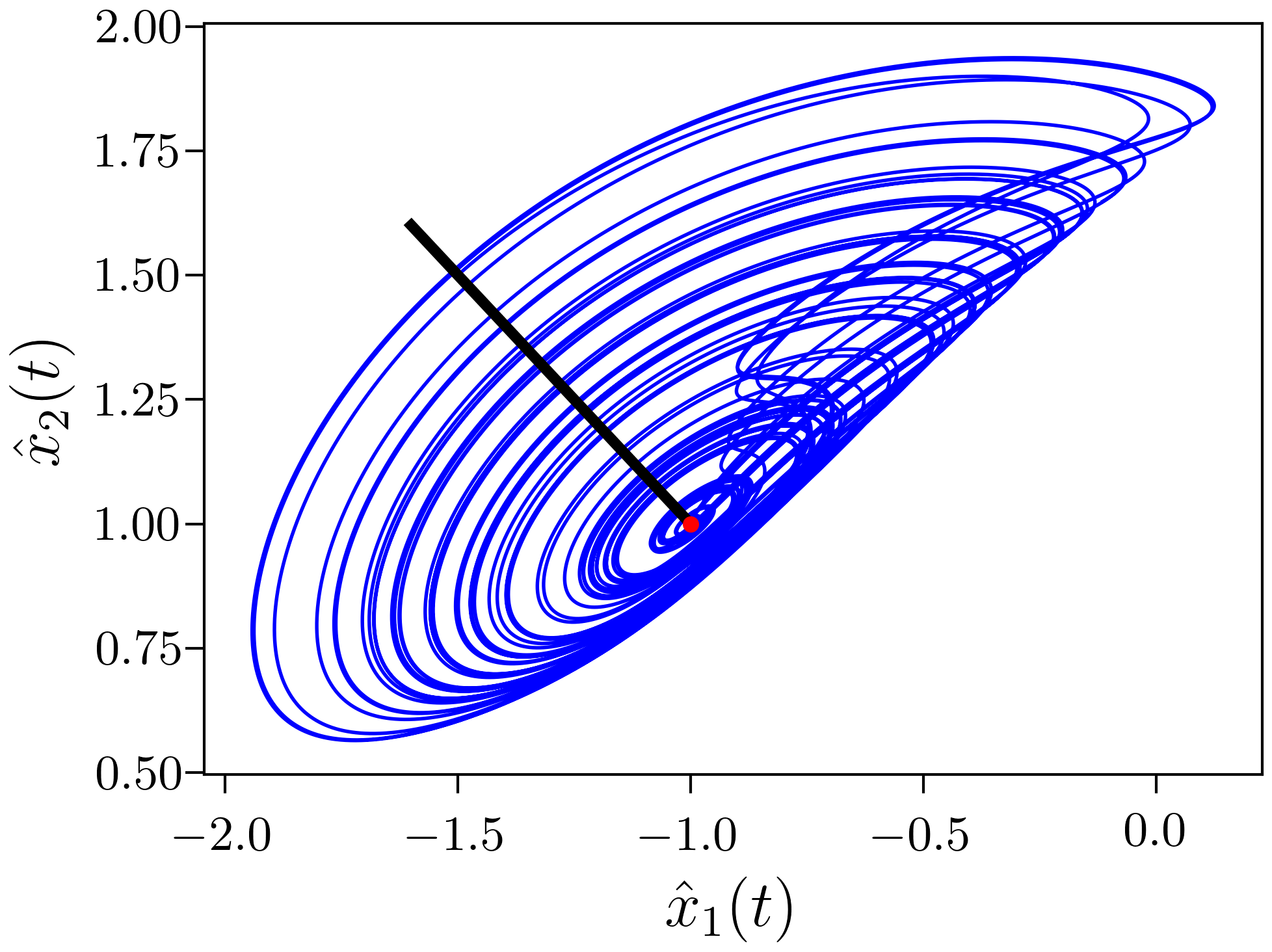} \ \includegraphics[width = 0.35\textwidth]{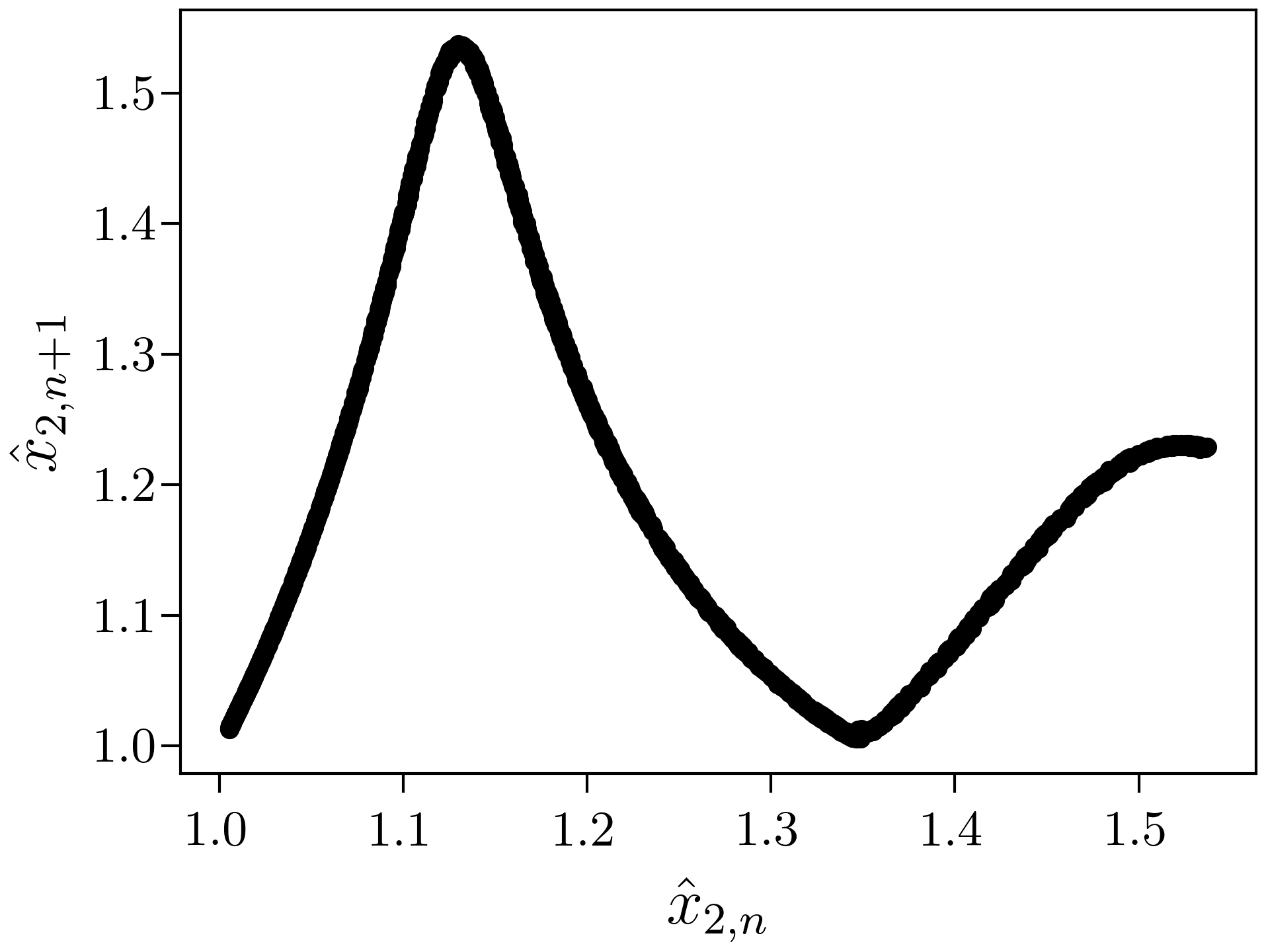}
\caption{Left: A typical chaotic trajectory of the Gissinger system \eqref{Gissinger} projected into the $(\hat x_1,\hat x_2)$-plane. In red is the equilibrium $(-1,1,0)$ and the black line represents a side view of the Poincar\'e section $\hat x_1 + \hat x_2 = 0$ and $\dot{\hat x}_1 + \dot{\hat x}_2 > 0$. Right: The first return map of the $\hat x_2$-variable on the attractor.}
\label{fig:Gissinger}
\end{figure} 

As an application of the dimensionality reduction capabilities of the neural network, we present another simple three-dimensional system known to exhibit chaotic dynamics. First presented by Gissinger~\cite{Gissinger}, we will consider the model
\begin{equation}\label{Gissinger} 
	\begin{split}
		\dot{x}_1 &= \mu x_1 - x_2x_3  \\
		\dot{x}_2 &= -\nu x_2 + x_1x_3 \\
		\dot{x}_3 &= \Gamma - x_3 + x_1x_2
	\end{split}
\end{equation}
Gissinger showed that by varying the parameters $(\mu,\nu,\Gamma) \in \mathbb{R}^3$ one would go from simple periodic motion to chaotic motion, most notably through a period-double cascade. The Poincar\'e map used in the original analysis is slightly tedious, therefore to simplify its definition we will shift and re-scale the $(x_1,x_2,x_3)$ variables appropriately. Applying the change of variable
\begin{equation}\label{ScaledGissinger}
	\begin{split}
		\hat{x}_1 &= x_1/\eta_1, \\
		\hat{x}_2 &= x_2/\eta_2 \\
		\hat{x}_3 &= (x_3 - \Gamma)/ \eta_3 \\
	\end{split}
\end{equation}
where $\eta_1 = \sqrt{\nu + \Gamma\sqrt{\frac \nu \mu}}, \eta_2 = \sqrt{\mu + \Gamma\sqrt{\frac \mu \nu}}$, and $\eta_3 = -\sqrt{\mu\nu} - \Gamma$, we move one equilibrium to the origin and scale the remaining two to be $(\pm 1,\mp 1,1)$. Then, Gissinger's Poincar\'e section becomes $\hat x_1 + \hat x_2 = 0$ and $\dot{\hat x}_1 + \dot{\hat x}_2 > 0$ in the hatted variables. \new{The reader should note that \eqref{ScaledGissinger} is an invertible change of variable and therefore is a conjugacy itself. This means that the dynamics in the Poincar\'e section of the hatted variables lie in one-to-one correspondence with the dynamics in Gissinger's section for \eqref{Gissinger} by simply applying a linear change of variable.} 

We will focus on the parameter values $\mu = 0.12$, $\nu = 0.1$, and $\Gamma = 0.85$, with the left panel of Figure~\ref{fig:Gissinger} presenting a typical chaotic trajectory at these parameter values. The section data is given in the variables $(\hat x_2,\hat x_3)$. The right panel of Figure~\ref{fig:Gissinger} presents the next iterate of $\hat x_2$ in the section against the current iterate. This data gives the appearance that the iterates of $\hat x_2$ in the section are almost entirely independent of $\hat x_3$. The next iterate of $\hat x_3$ in the section against the current iterate looks similar but is not provided here for brevity. We have calculated the LEs to two decimal places to be $\lambda_1 = 0.07$, $\lambda_2 = 0$, and $\lambda_3 = -1.05$. This gives a Kaplan--York dimension of $2.07$, and from the discussion above we would expect the dynamics in the section to be conjugate to a one-dimensional mapping.

We employ the neural network to discover a conjugacy of the Gissinger Poincar\'e section data with a one-dimensional homogeneous cubic map.
Using the specifications described in the appendix we arrive at a conjugacy with the mapping
\begin{equation}\label{GissingerMap}
	g(y) = 8.5344y -18.2999y^2 + 9.8172y^3
\end{equation} 
The election to conjugate to a cubic map is motivated by the appearance of the training data in Figure~\ref{fig:Gissinger}. In much the same way that the R\"ossler system gives the appearance that its section data is governed by a quadratic map, the Gissinger system gives the appearance that its section is governed by a cubic map. The difference is that the section data is governed by two variables, $(\hat x_2,\hat x_3)$, and so attempting to discover a mapping governing these iterates would either require a two-dimensional mapping or losing information about the potential slight influence the variables have upon each other. Hence, the advent of using the network to discover conjugacies is that it can automate the task of dimensional reduction while also providing a nonlinear change of variable to fit to a simple model. Moreover, the simplicity of the mapping \eqref{GissingerMap} lends itself well to extracting UPOs and mapping them back to UPOs in the Gissinger system \eqref{Gissinger}. To provide evidence for this fact, the accompanying supplementary material contains a MATLAB script to stabilize UPOs that intersect the Poincar\'e section once, twice, three times, and four times.  

\new{Prior to concluding this subsection we comment that additional explorations of the Gissinger system \eqref{Gissinger} were undertaken to verify the consistency of our results. This includes working with the original variables $(x_1,x_2,x_3)$ instead of introducing the change of variable \eqref{ScaledGissinger}. In this case the coefficients of the discovered conjugate cubic mapping were within $4.5\%$ of those of \eqref{GissingerMap}. One method that we chose to compare the results between the discovered conjugacy trained on the Poincar\'e section data with and without the change of variable \eqref{ScaledGissinger} is to compare the location of the fixed point (corresponding to the lowest period UPO) that they give in the original Gissinger model \eqref{Gissinger}. To do this we solve $g(y) = y$ for the respective conjugacy mapping and then map the result back to the Gissinger system using the associated inverse of the conjugacy, $\hv^{-1}$. The resulting location from the network trained on the original Gissinger variables is within $0.3\%$ of the location from the neural network trained on the hatted variables \eqref{ScaledGissinger}. Similar conjugacies to mappings nearly identical to \eqref{GissingerMap} have also been established when using different Poincar\'e sections.}  

\section{Infinite Dimensional Systems}\label{sec:Infinite} 

Having now demonstrated the utility of the neural network on low-dimensional chaotic systems, we turn now to some more complex examples. In particular, we will focus on the Kuramoto--Sivashinsky PDE and the Mackey--Glass equation. In the case of the former we demonstrate how we can use high-dimensional approximations of the system coming from Galerkin projections to identify the low-dimensional dynamics on the attractor. From these low-dimensional dynamics we can then extract the intersections of the UPOs with the Poincar\'e section and implement the methods of Subsection~\ref{subsec:UPO} to obtain the UPOs in the continuous-time system. In the case of the Mackey--Glass equation we show how two different Poincar\'e maps are conjugate to the same simple quadratic mapping and therefore can be used to provide evidence for a long-standing conjecture surrounding the infinite-dimensional system.

\subsection{Kuramoto--Sivashinsky Equation}\label{sec:KS} 

\begin{figure} 
\center
\includegraphics[width = 0.3\textwidth]{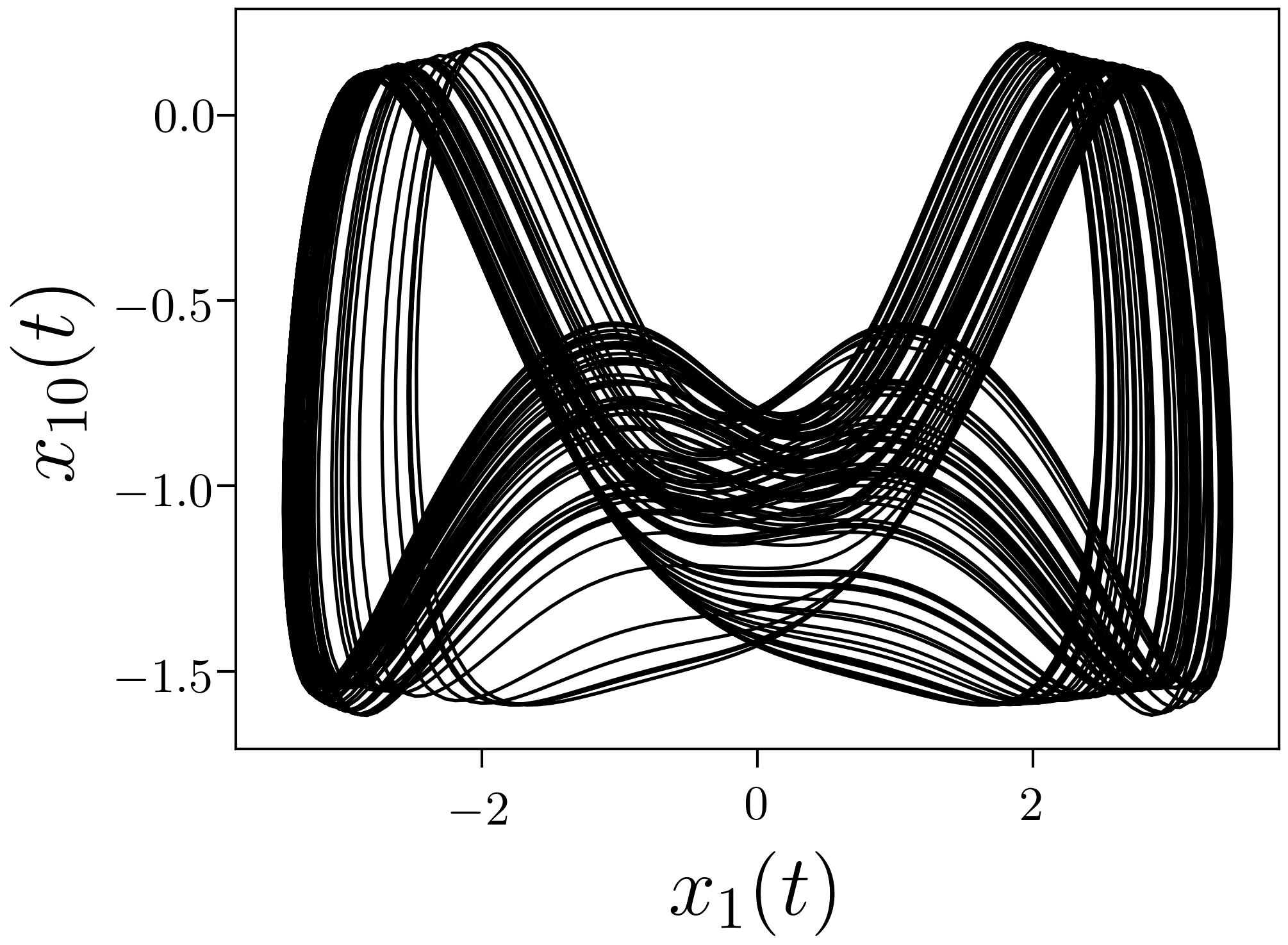} \ \includegraphics[width = 0.3\textwidth]{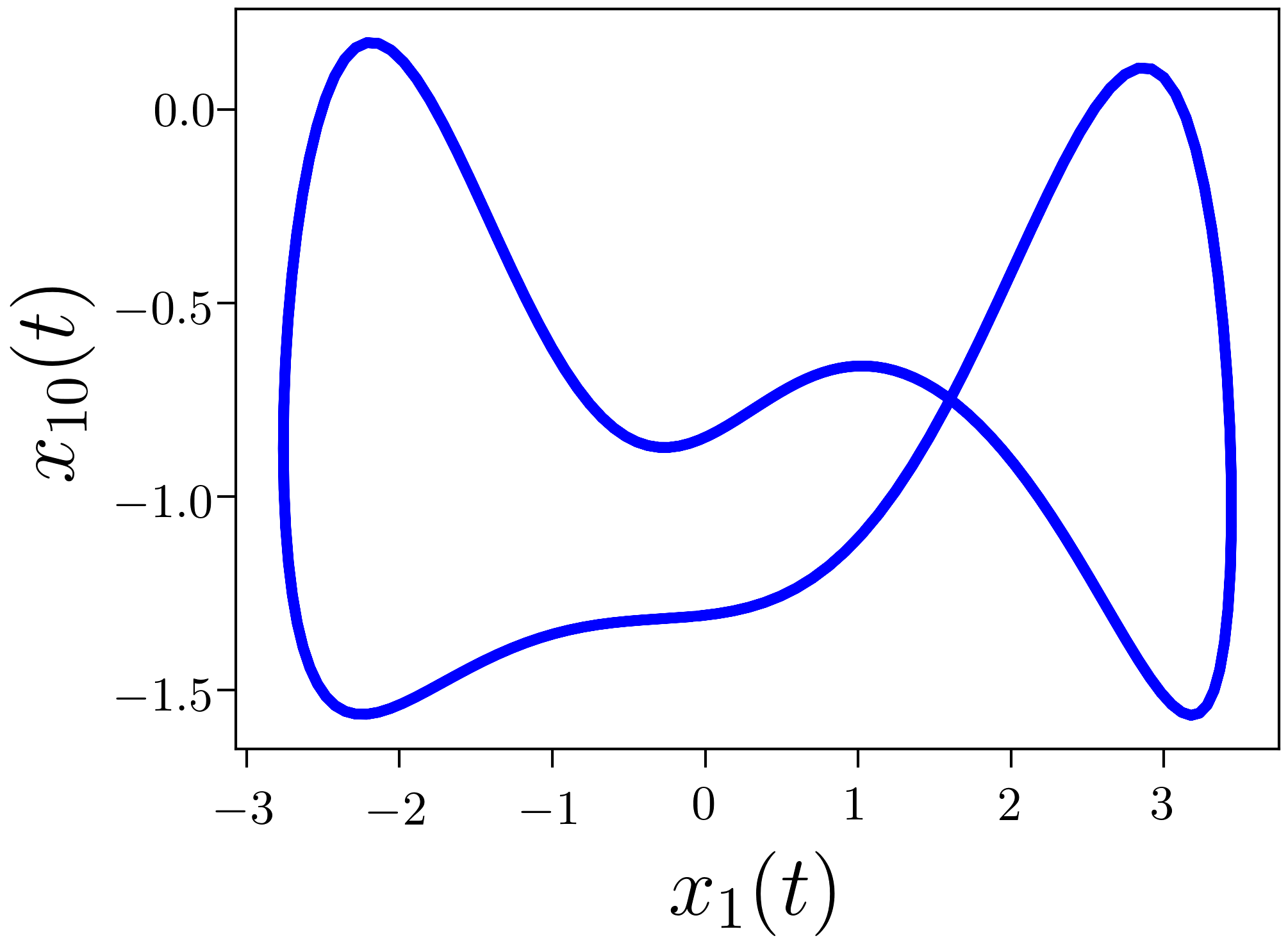} \ \includegraphics[width = 0.3\textwidth]{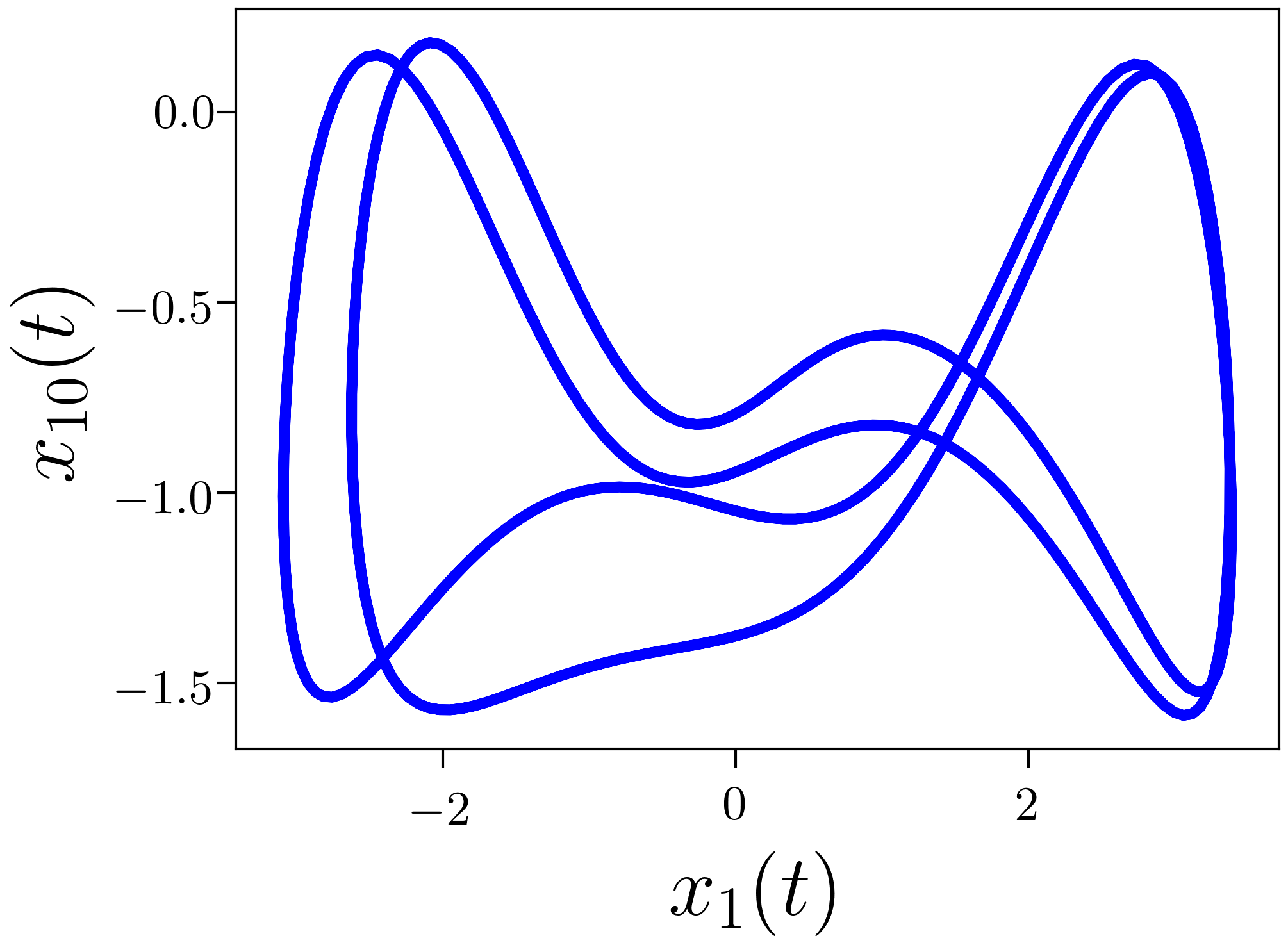} \ \includegraphics[width = 0.3\textwidth]{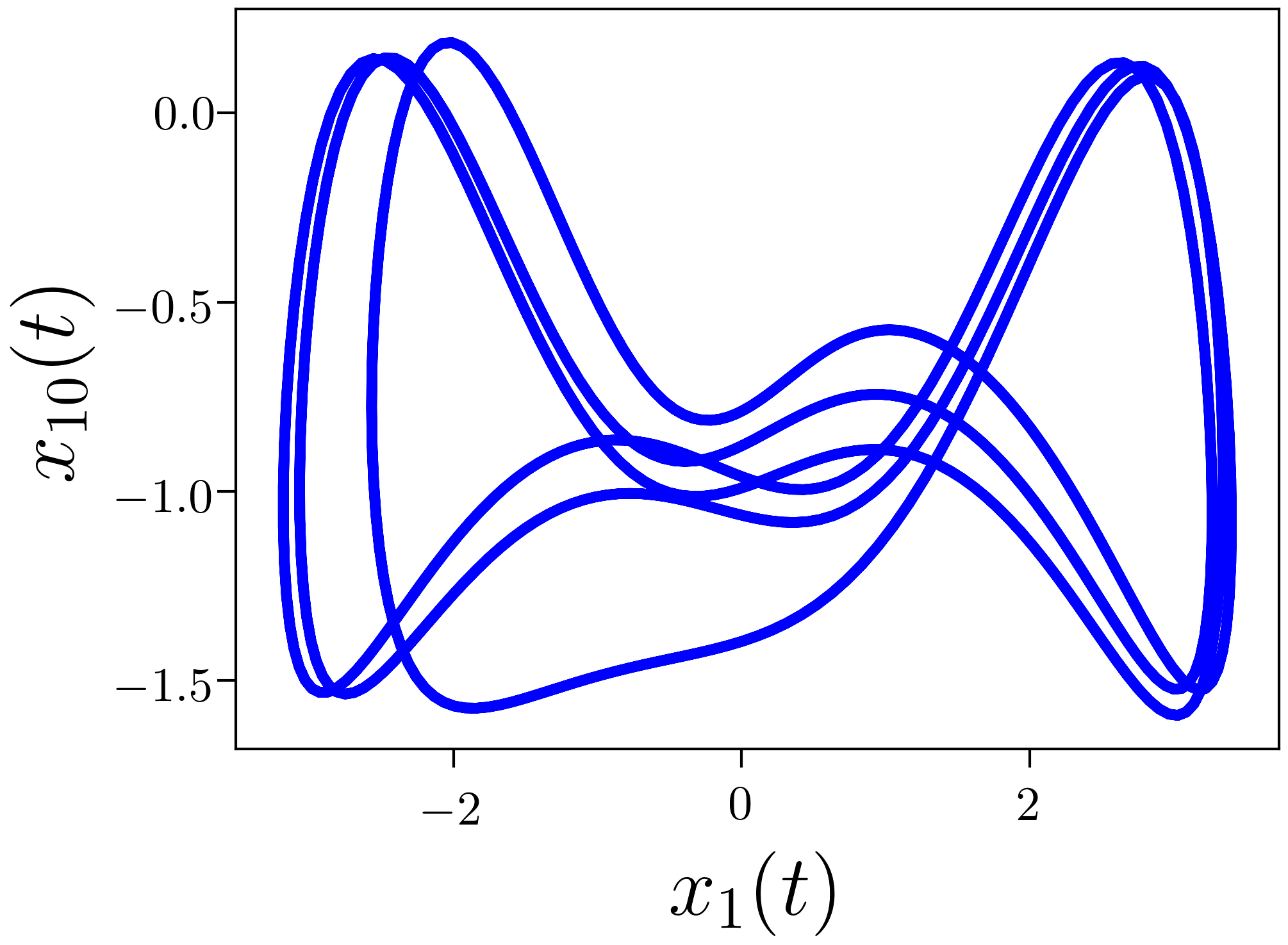} \ \includegraphics[width = 0.3\textwidth]{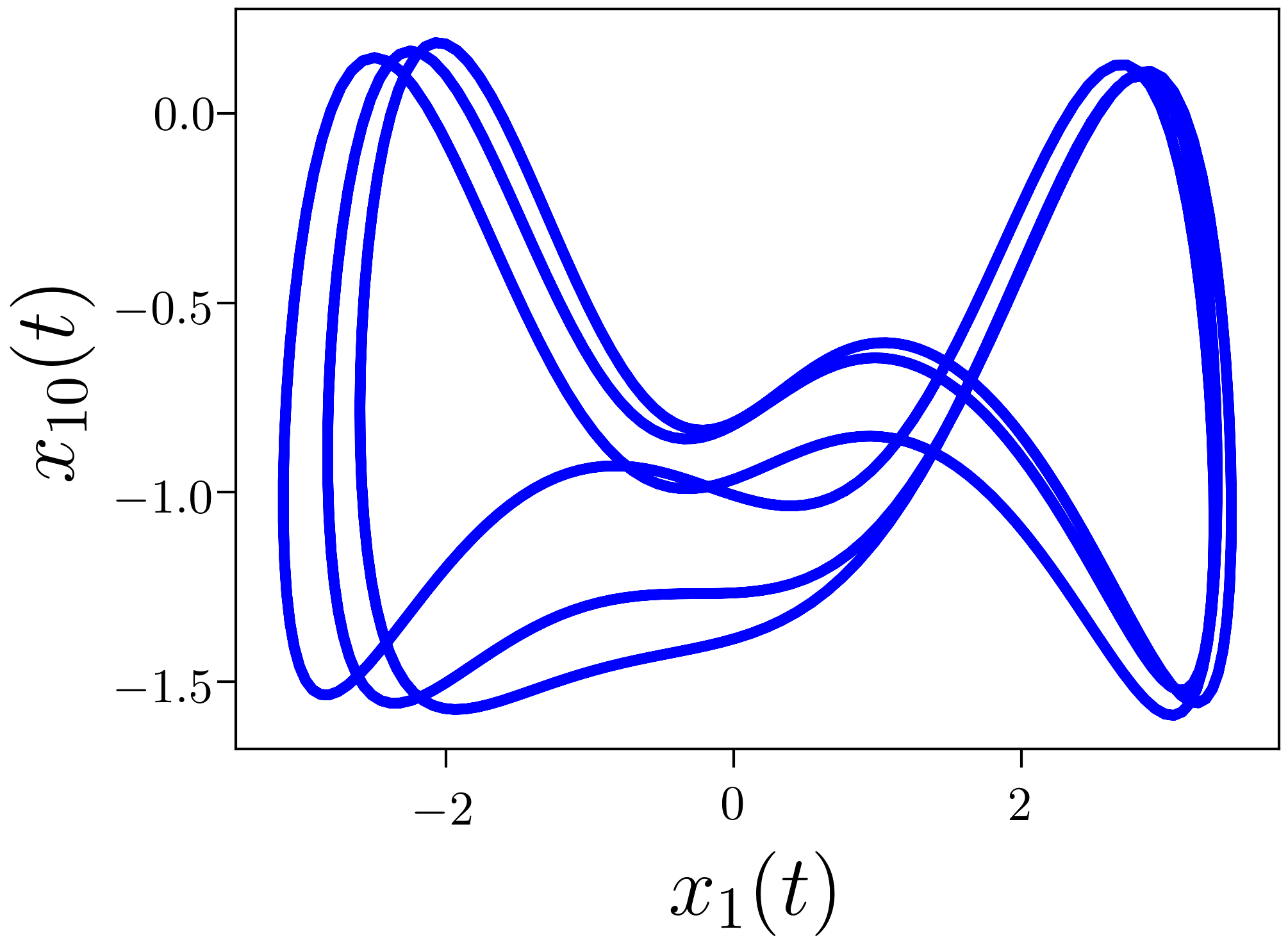}  \ \includegraphics[width = 0.3\textwidth]{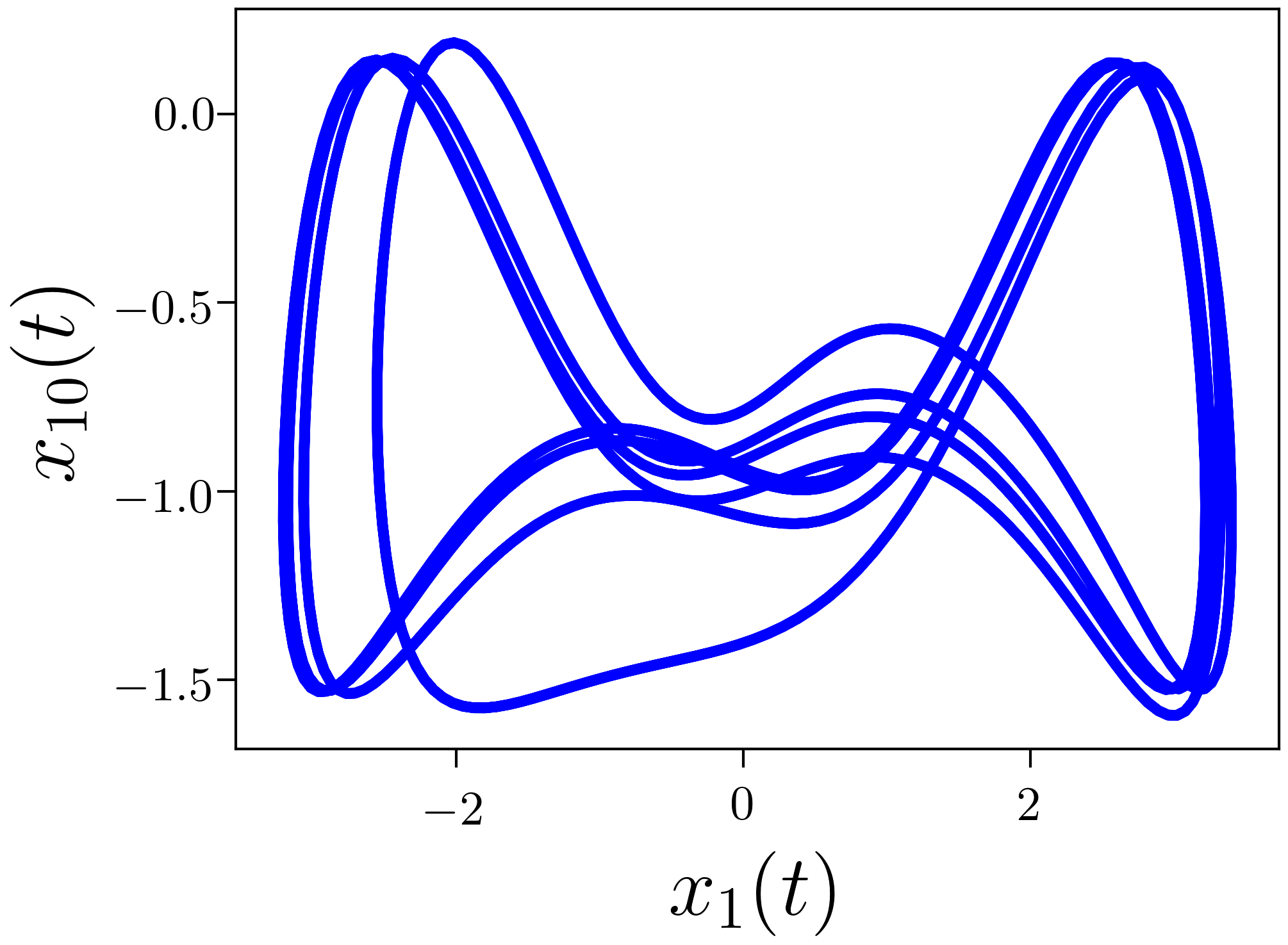}
\caption{The chaotic attractor (black) of \eqref{Kuramoto} with $\nu = 0.0298$ projected into the $(x_1,x_{10})$-plane along with identified UPOs (blue). Top row: R and LR orbits. Bottom row: LLR, LRR, and LLLR.}
\label{fig:KSE}
\end{figure}

Here we will consider the Kuramoto--Sivashinsky equation (KSE), given by the PDE
\begin{equation}\label{Kuramoto}
	u_t + \nu u_{\xi\xi\xi\xi} + u_{\xi\xi} + uu_\xi = 0,
\end{equation}
where $u = u(\xi,t)$ is a function of space $\xi \in [-\pi,\pi]$ and time $t \geq 0$. Following~\cite{Kuramoto,Kuramoto2}, we restrict ourselves to the flow-invariant set of odd periodic functions $u(\xi,t)$. We perform the Galerkin/Fourier projection onto $N \geq 1$ spatial modes with time varying coefficients 
\begin{equation}
	u(\xi,t) = \sum_{k = 1}^N x_{k}(t)\sin(k\xi).
\end{equation} 
This leads to an $N$-dimensional coupled system of differential equations for the time-dependent coefficients, given by
\begin{equation}\label{KSE2}
	\dot x_k = k^2(1 - \nu k^2)x_k + \frac{k}{2}\sum_{i = 1}^{N-k}x_ix_{k+i} - \frac{k}{4}\sum_{j = 1}^{k}x_j x_{k-j}, \quad k = 1,\dots, N.
\end{equation}
In what follows we will fix $N = 14$ and generate training data with initial conditions $x_k(0)$ drawn from the uniform distribution on $[0,0.1]$. The top left image in Figure~\ref{fig:KSE} presents a typical chaotic trajectory of \eqref{KSE2} with $\nu = 0.0298$ projected into the $(x_1,x_{10})$-plane. 

Previous numerical investigations of the KSE using the system \eqref{KSE2} have shown that a period-doubling sequence into chaos takes places as $\nu$ descends towards ${\sim}0.0297$. Furthermore, the works~\cite{Kuramoto,Kuramoto2} calculate that this period-doubling sequence obeys the Feigenbaum scaling - a universal feature of one-dimensional maps with a quadratic extremum. Notably, these works draw frequent comparisons to the logistic map \eqref{Logistic}, but lack an analytical backing for such a correspondence. Here we will follow~\cite{Kuramoto3} and define a Poincar\'e section by the half-plane $x_1 = 0$ and $\dot x_1 > 0$ to produce a 13-dimensional collection of discrete section data. With the neural network architecture described above, we are now in a position to obtain a conjugacy between this 13-dimensional section data and the logistic map.

To illustrate our results we will take $\nu = 0.0298$, which is near the end of the period-doubling cascade and represents a region of chaotic dynamics in the system. The three largest LEs are $\lambda_1 = 0.74$, $\lambda_2 = 0$, and $\lambda_3 = -5.92$, thus giving a Kaplan--Yorke dimension of $2.13$. This provides further justification for attempting to obtain a conjugacy with a one-dimensional mapping. We find a conjugacy with the quadratic mapping 
\begin{equation}\label{KSMap}
	g(y) = 3.9653y -3.9153y^2,
\end{equation}       
and comment that applications of the method with the addition of cubic and/or quartic terms to the mapping show little improvement over the obtained quadratic model. From the discussion on the R\"ossler attractor, the discovered latent map \eqref{KSMap} generates equivalent dynamics to the logistic map \eqref{Logistic} with $r = 3.9653$. As expected from the numerical results, we find that the chosen value $\nu = 0.0298$ is near the culmination of the period-doubling cascade and firmly in a chaotic parameter region.

The utility of the discovered mapping \eqref{KSMap} is that it provides a low-dimensional mapping that can be used to understanding the chaotic attractor in the infinite-dimensional system \eqref{Kuramoto}. To demonstrate this utility, we extract periodic orbits from \eqref{KSMap} and relate them to the 13-dimensional Poincar\'e section via the neural-network conjugacy. These section points are then used to seed initial guesses to produce UPOs which form the skeleton of the chaotic attractor in the same way that was done for the Lorenz system above. Furthermore, we search for the UPOs in system \eqref{KSE2} with $N = 28$, double the number of variables used for the training data, to improve accuracy. The additional 14 variables are simply initialized as zero functions, aiding in showing that our 14 mode truncation indeed captures much of the dynamical structure of the chaotic system. We label the identified UPOs using the symbols L and R, meaning that the iterate of the map \eqref{KSMap} is either to the left or the right of its global maximum, respectively. Our results are summarized in table~\ref{tbl:KSE}, where we emphasize how few iterates of the solver are required to produce a UPO solution to \eqref{KSE2}. We provide data on all periodic orbits with sequence length up to five and comment that the sequences LLLLR and LLLRR are notably absent since such periodic orbits are not present in \eqref{KSMap}. Figure~\ref{fig:KSE} presents visual demonstrations of some of these identified UPOs projected onto the $(x_1,x_{10})$-plane. 

\begin{table}
\centering
\renewcommand{\arraystretch}{1}
\begin{tabular}{ |c||c|c|c| }
\hline
Symbolic Sequence & Period & fsolve Iterates & fsolve Time (seconds)  \\
\hline
L & 0.8630 & 7 & 4.9900 \\
R & 0.8567 & 6 & 6.0526 \\
LR & 1.7356 & 9 & 37.5210 \\
LLR & 2.6242 & 10 & 96.2644 \\
LRR & 2.6031 & 8 & 61.9304 \\
LLLR & 3.5167 & 10 & 167.6217  \\
LLRR & 3.5166 & 9 & 124.7951 \\
LRRR & 3.4623 & 10 & 159.5762 \\
LLRLR & 4.3651 & 11 & 441.0552 \\
LLRRR & 4.3539 & 12 & 306.9639 \\
LRRLR & 4.3424 & 11 & 236.0808 \\
LRRRR & 4.3247 & 9 & 210.8622 \\
\hline
\end{tabular}
\caption{Symbolic sequences of periodic orbits in the discovered conjugate mapping are used to obtain periodic orbits in the KSE \eqref{Kuramoto}. All symbolic sequences that are present in the discovered mapping up to length 5 are given, along with the same information provided in table~\ref{tbl:Lorenz}.}
\label{tbl:KSE}
\end{table}

As $\nu > 0$ is lowered beyond the culmination of the period-doubling cascade the system \eqref{Kuramoto} exhibits a variety of exotic behaviour that is not fully understood~\cite{Kuramoto,Kuramoto2}. In particular, the chaotic attractor becomes increasingly complicated, indicated by larger LEs and an increasing Kaplan--Yorke dimension. The further the Kaplan--Yorke dimension of the attractor is above 2, the more difficulty we have in understanding its topology. Therefore, the choice for the form of the latent map $g$ becomes less clear since we cannot easily relate the Poincar\'e section dynamics to well-studied universal maps.  

\begin{figure} 
\center
\includegraphics[width = 0.3\textwidth]{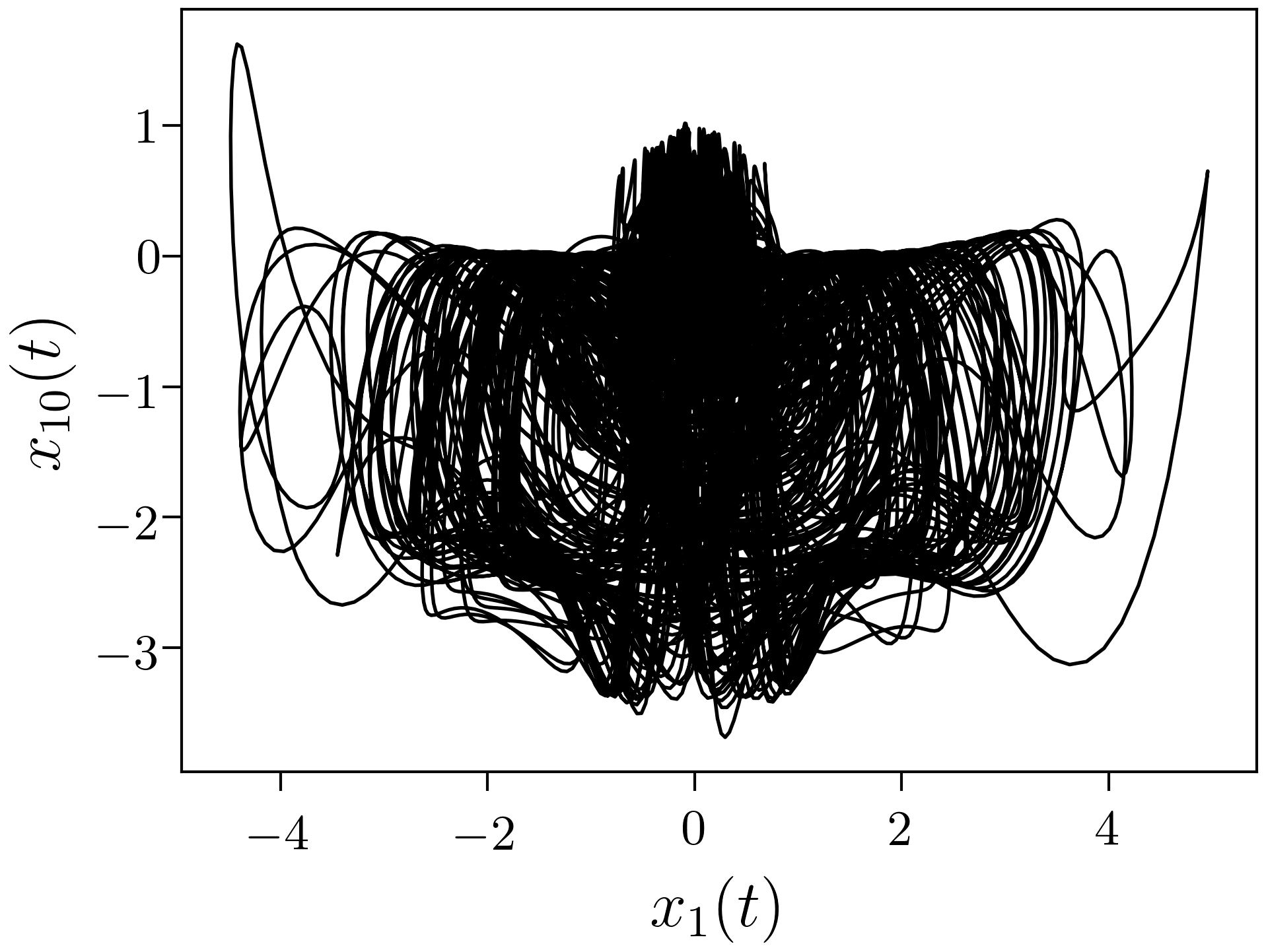} \ \includegraphics[width = 0.3\textwidth]{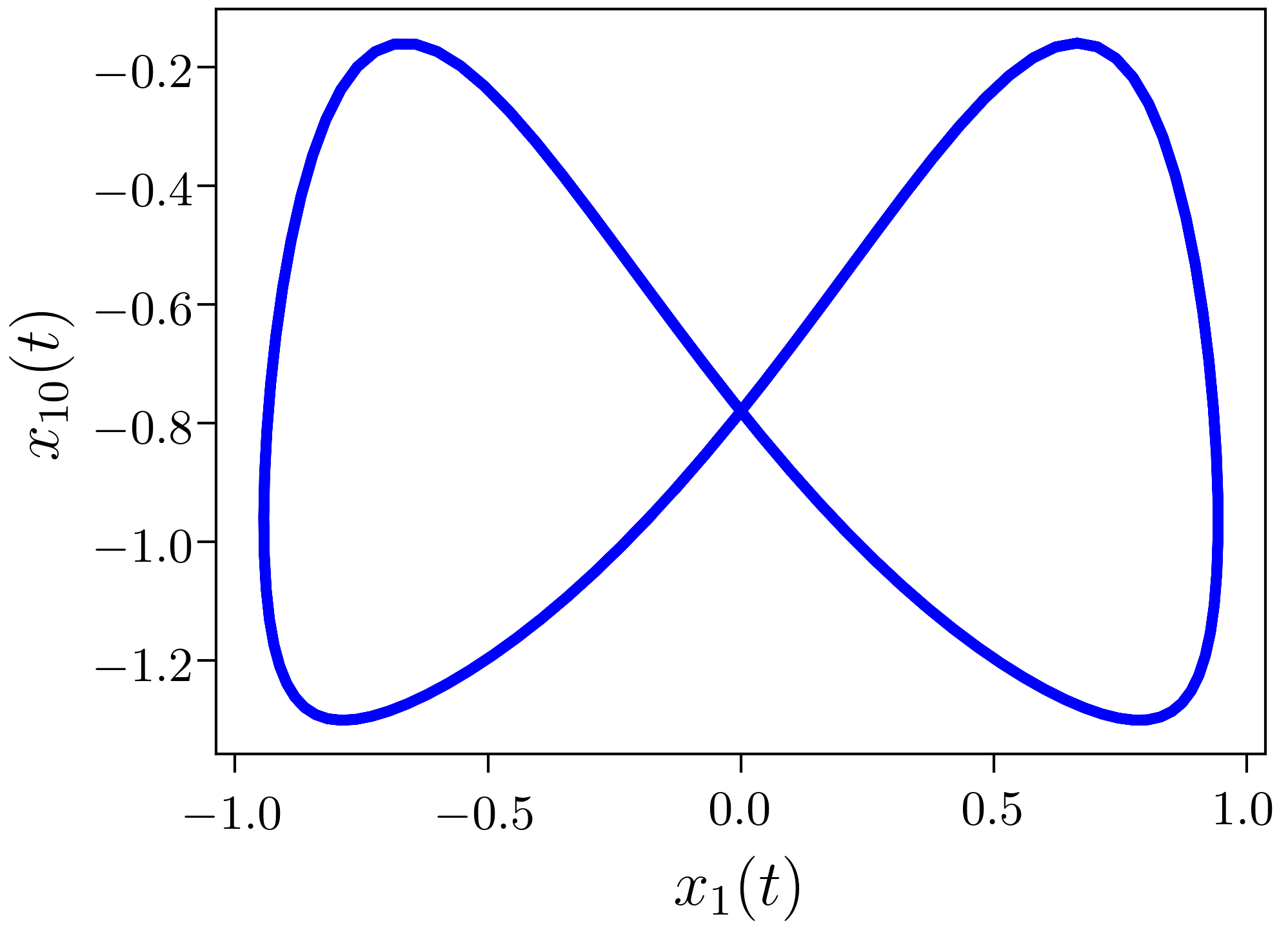} \\ \includegraphics[width = 0.3\textwidth]{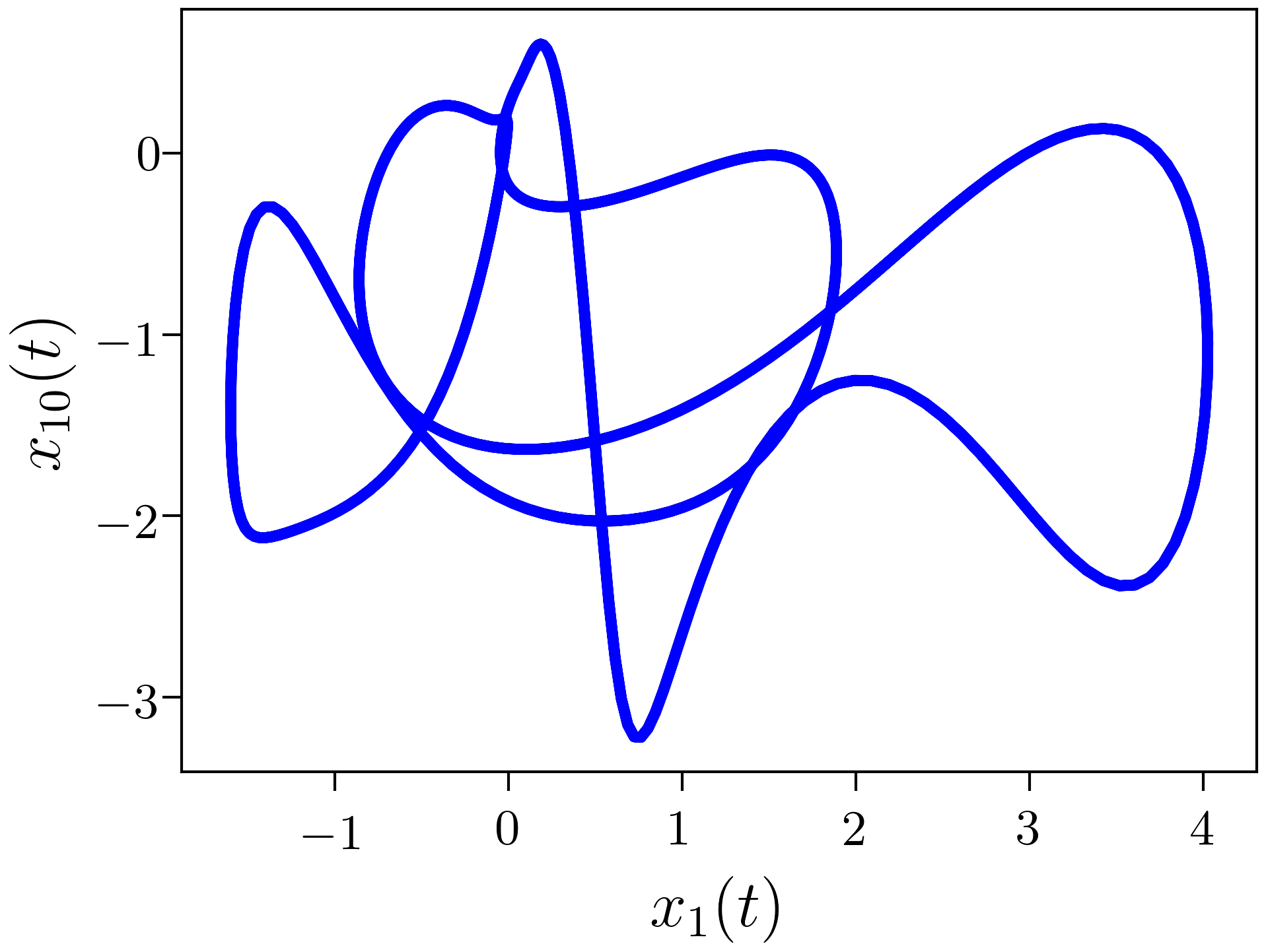} \ \includegraphics[width = 0.3\textwidth]{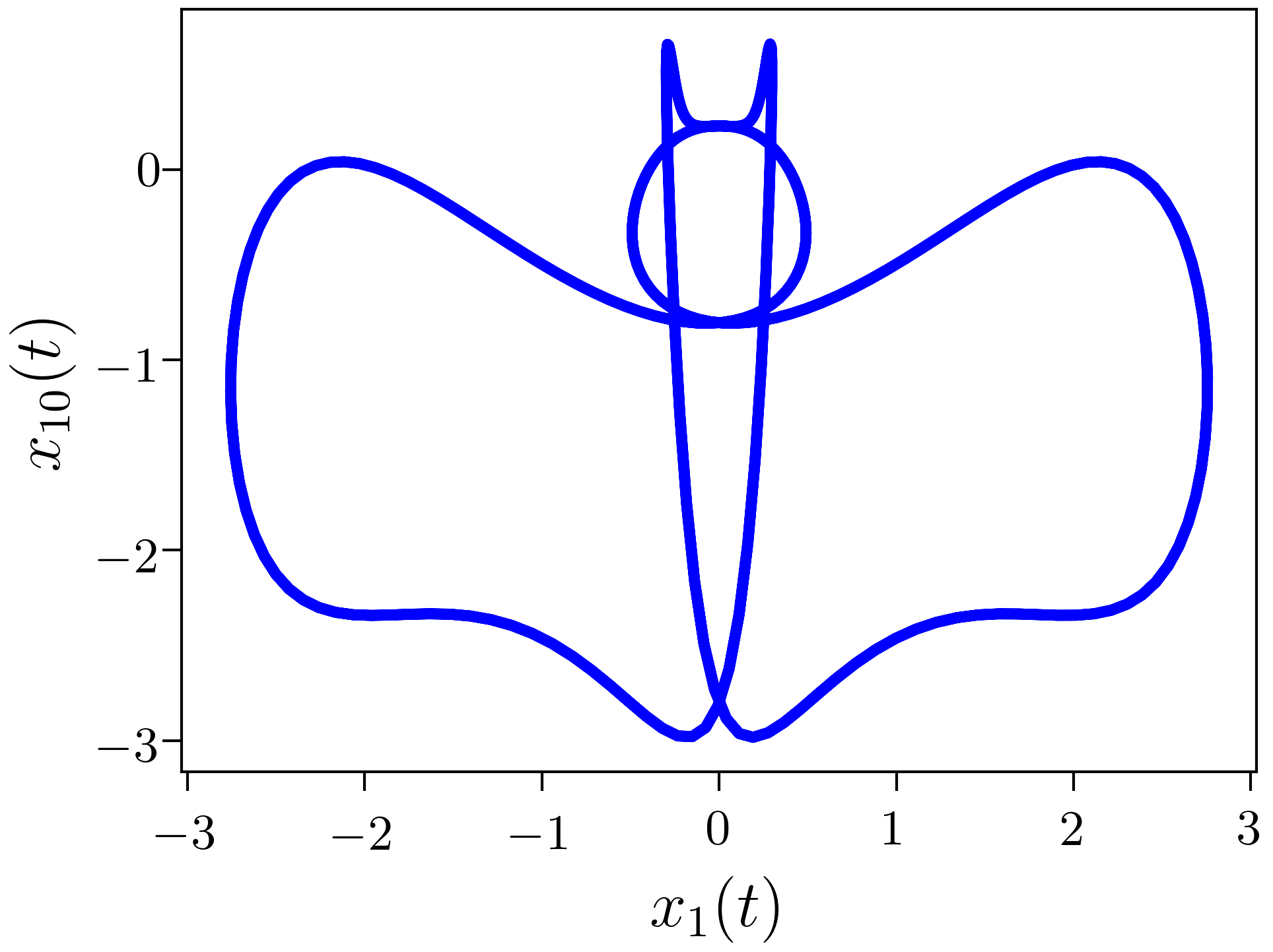} 
\caption{The chaotic attractor (black) of \eqref{Kuramoto} with $\nu = 0.0210$ projected into the $(x_1,x_{10})$-plane along with identified UPOs (blue). Top right: A UPO that intersects the Poincar\'e section exactly once. Bottom: UPOs that intersect the Poincar\'e section exactly twice and three times, respectively.}
\label{fig:KSE2}
\end{figure}

Taking $N = 14$ and $\nu = 0.021$, the first four LEs of the KSE \eqref{KSE2} are given by
\begin{equation}
	\lambda_1 = 1.75, \quad \lambda_2 = 0, \quad \lambda_3 = -1.65,\quad \lambda_4 = -4.85,
\end{equation}
giving a Kaplan--Yorke dimension of $3.02$. In the top left of Figure~\ref{fig:KSE2} we provide a typical chaotic trajectory of the system at $\nu = 0.021$. Based on the size of the Kaplan--Yorke dimension, we seek a latent mapping that is two-dimensional and conjugate to the 13-dimensional Poincar\'e section data. The network is able to identify a conjugacy with the quadratic function
\begin{equation}
	\gv(y_1,y_2) = \begin{bmatrix}
	0.1009 + 1.5589y_1 + 0.5601y_2 -0.2534y_1^2 -2.4928y_1y_2 + 0.4457y_2^2 \\
	0.6545 + 0.4098y_1 -1.2802y_2 + 0.2022y_1^2 + 0.6164y_1y_2 -0.3460y_2^2   
	\end{bmatrix}.
\end{equation}
We again demonstrate the utility of this mapping by using it to identify UPOs in the chaotic flow of \eqref{KSE2}. We use a 28-mode truncation to identify UPOs with the final 14 modes initialized identically to zero. Figure~\ref{fig:KSE2} shows plots of three UPOs projected into the $(x_1,x_{10})$-plane. The presented UPOs intersect the Poincar\'e section at exactly one point, two points, or three points. 

\subsection{Mackey--Glass Equation}\label{sec:Mackey} 

\begin{figure} 
\center
\includegraphics[width = 0.35\textwidth]{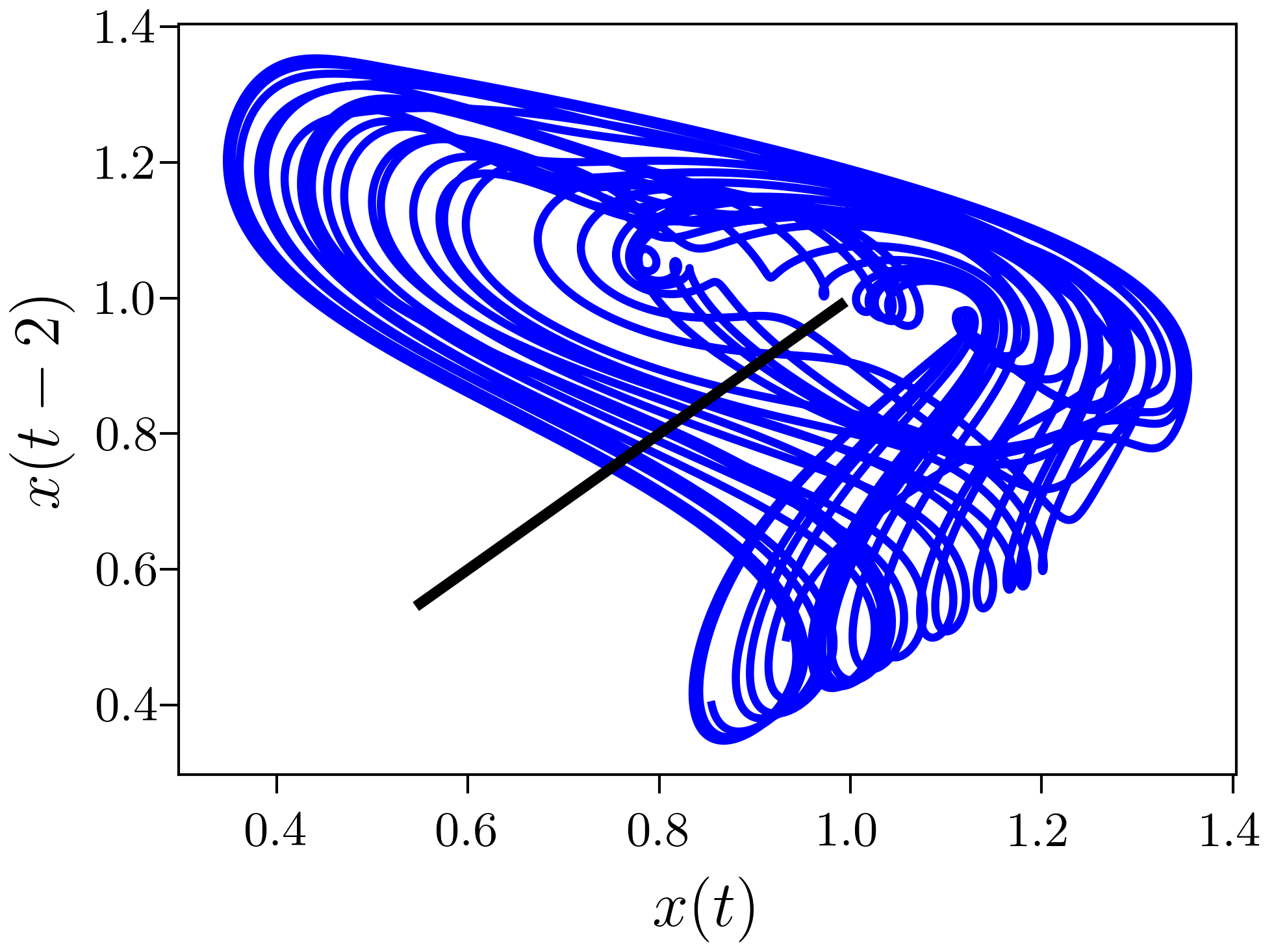} \ \includegraphics[width = 0.35\textwidth]{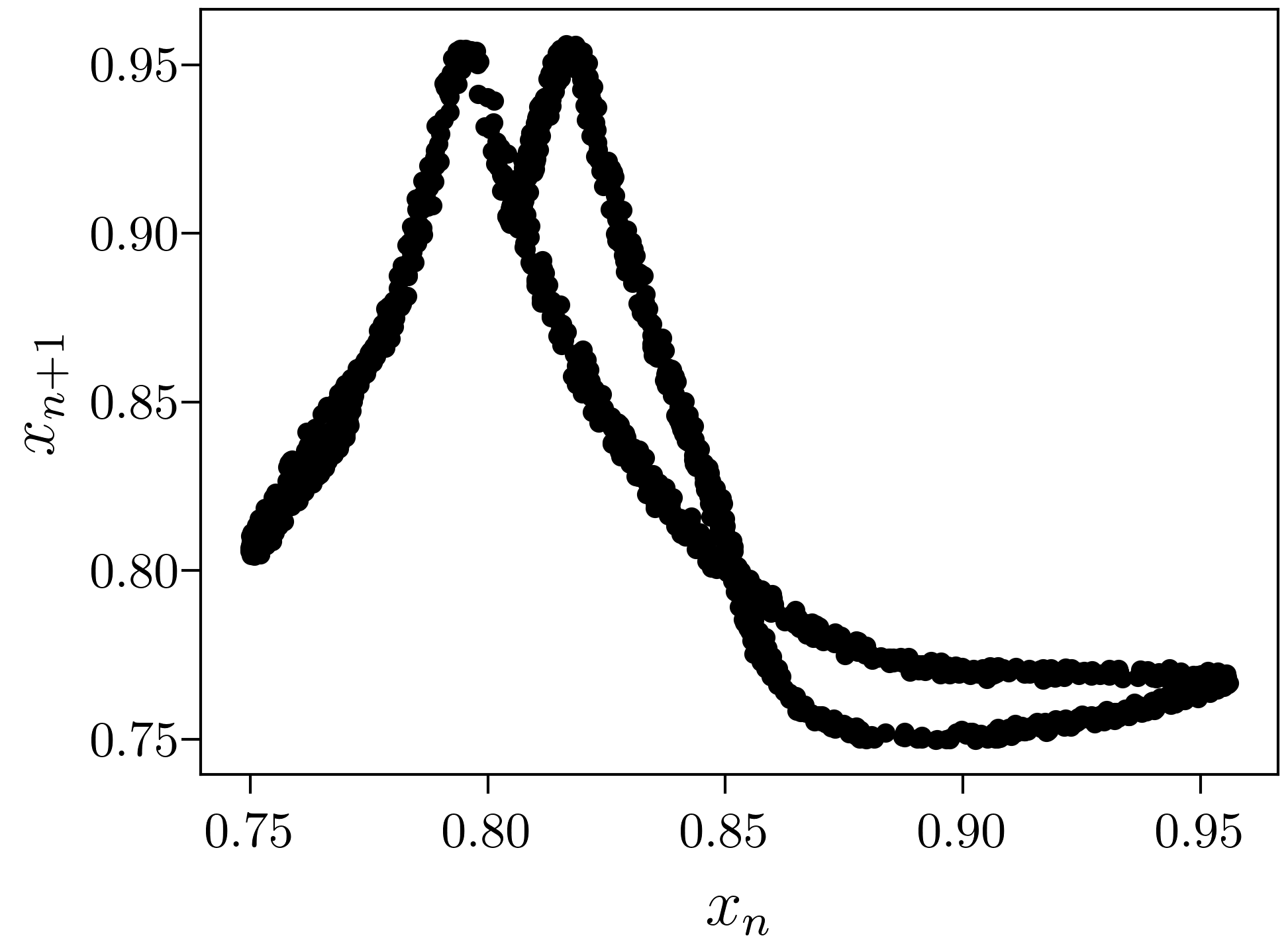} \\
\includegraphics[width = 0.35\textwidth]{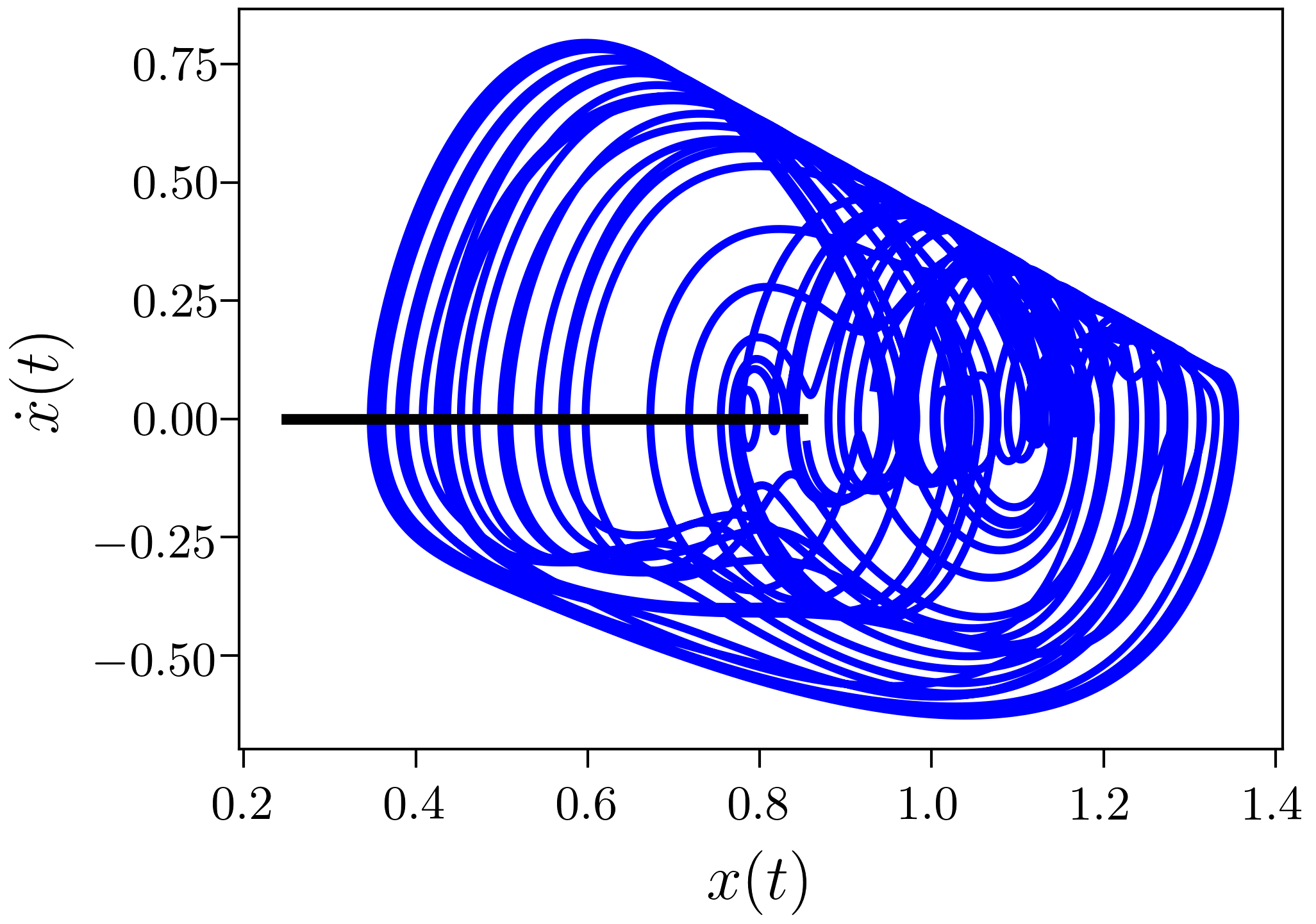} \ \includegraphics[width = 0.35\textwidth]{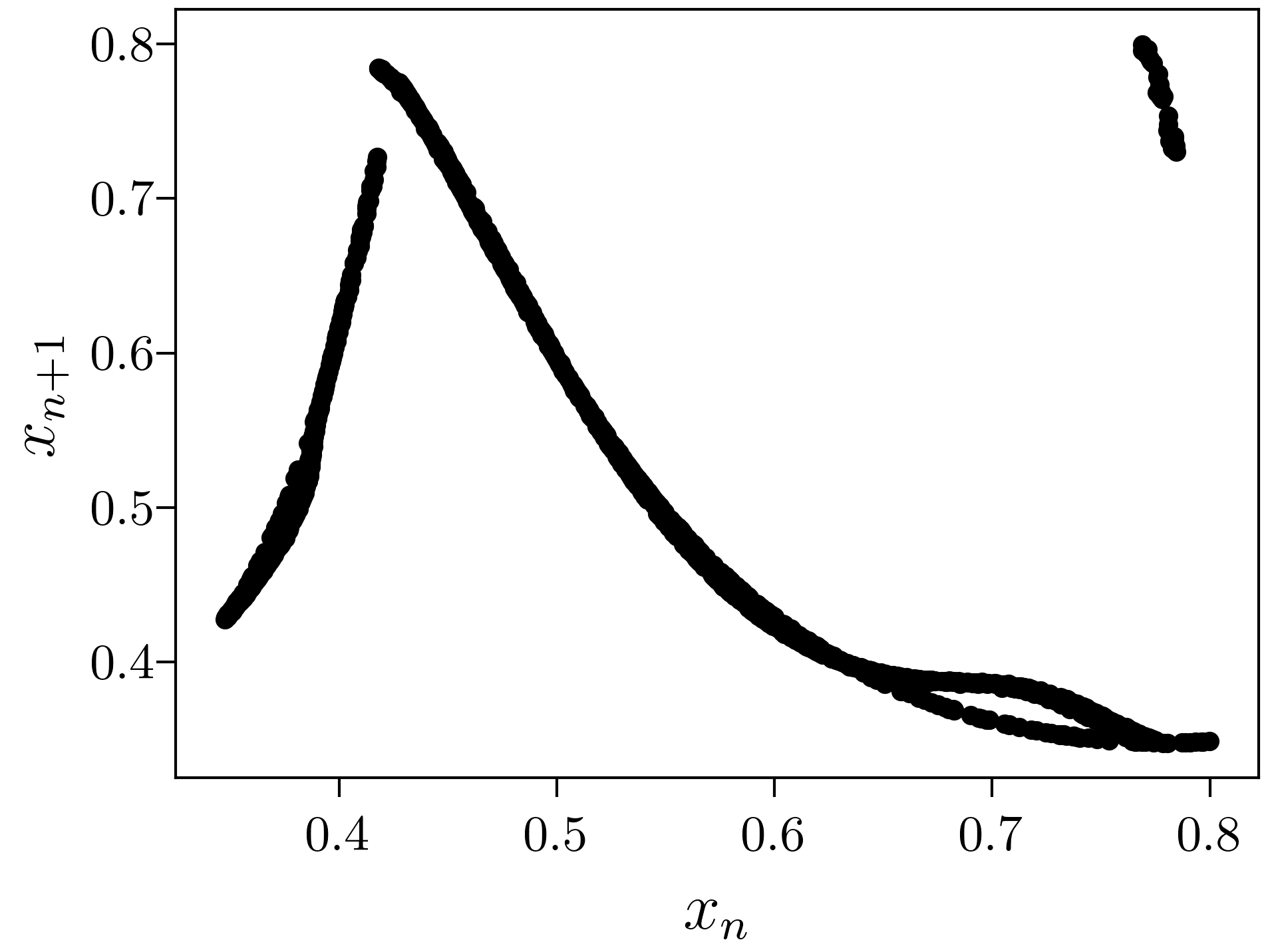}
\caption{Left: A typical chaotic trajectory of the Mackey--Glass delay equation \eqref{Mackey} with $p = 9.65$ projected into the (top) $(x(t),x(t-2))$-plane and into the (bottom) $(x(t),\dot{x}(t))$-plane. The black line represents the Poincar\'e section (top) $x(t) = x(t-2)$ and (bottom) $\dot{x}(t) = 0$. Right: The first return maps of the $x$-variable on the attractor for the respective Poincar\'e sections shown on the left.}
\label{fig:Mackey}
\end{figure} 

Let us now consider the Mackey-Glass delay differential equation~\cite{Mackey}
\begin{equation}\label{Mackey}
	\dot{x} = \frac{2x_\tau}{1 + x_\tau^p} - x,
\end{equation}
where $x_\tau = x(t - \tau)$ is the value of the function $\tau > 0$ time units in the past. We will fix $\tau = 2$ throughout this investigation. Then, increasing $p$ from approximately $7$ results in a period-doubling cascade into chaos~\cite{Mackey}. On the left of Figure~\ref{fig:Mackey} we present typical chaotic trajectories with $p = 9.65$ in both the $(x(t),x(t-2))$-plane and the $(x(t),\dot{x}(t))$-plane. This sequence of period-doubling bifurcations has led to conjectures on whether there is a correspondence between the well-known sequence of bifurcations leading to chaos in quadratic maps, such as the logistic equation \eqref{Logistic}, and the infinite-dimensional Mackey--Glass equation \eqref{Mackey}. This remains an open problem that we will provide numerically-assisted evidence for.

There are two different Poincar\'e sections that are used to analyze the infinite-dimensional dynamics of \eqref{Mackey}. They are given by~\cite{Scholarpedia} 
\begin{equation}\label{MGPsec}
	\begin{split}
	 (1)& \quad x(t) = x(t-2),\quad s.t.\quad \dot{x}(t) - \dot{x}(t - 2)\geq 0,\ x(t) < 0.96 \\
	 (2)& \quad \dot{x}(t) = 0,\quad s.t.\quad \ddot{x}(t) > 0,\ x(t) < 0.8.
	 \end{split}
\end{equation}
These sections are illustrated in Figure~\ref{fig:Mackey}, along with their respective first return maps. To illustrate how correspondences between these two distinct Poincar\'e maps and the logistic map \eqref{Logistic} can be found using the network we will fix $n = 9.65$ and show that both maps correspond to the same value of logistic map parameter, denoted $r$. Indeed, at $n = 9.65$ the network produces conjugacies with the quadratic maps
\begin{equation}
	\begin{split}
		(1)&\quad g_1(y) = 3.8390y - 3.9016y^2 \\
		(2)&\quad g_2(y) = 3.8321y -3.8976y^2,
	\end{split}
\end{equation}
where (1) and (2) denote the choices of Poincar\'e maps in \eqref{MGPsec}.  From the discussion on the R\"ossler system, the maps are equivalent to the logistic map with $r = 3.8390$ and $r = 3.8321$, respectively.

\section{Discussion}\label{sec:Discussion} 

Through a number of illustrative examples, we have seen how the proposed deep learning architecture can be used to understand and classify chaotic behaviour by learning coordinates and dynamics that can be leveraged to understand conjugacies. Indeed, these examples have demonstrated that by using the network to learn conjugate mappings we can improve forecasting of chaotic systems, achieve a dimensionality reduction from the observed variables, and obtain actionable and interpretable mappings. These aspects of the method were highlighted with three low-dimensional chaotic systems: the R\"ossler, Lorenz, and Gissinger systems. Then, with an improved intuition and understanding of the method, we applied it to two infinite-dimensional systems: the Kuramoto--Sivashinsky PDE and the Mackey--Glass delay-differential equation. It was evident from the analysis of the Kuramoto--Sivashinsky PDE that this method can go a long way towards better understanding spatio-temporal chaos. It therefore becomes appealing to apply these methods to fluid systems, where the dimensionality reduction can particularly help to understand low-dimensional turbulent behaviour in extremely high-dimensional systems. Results in this direction will be reported in a follow-up investigation to this one.

\new{We emphasize that our methods are comparable to other data-driven discovery algorithms in that the goal is to obtain an explicit dynamical system by expanding an unknown function in a library of candidate functions. Hence, the discovery process amounts to tuning the coefficients of the linear span for which the discovered mapping is assumed to belong to. The advantage of using the autoencoder structure proposed in our work is that finding the change of variable, $\hv$, and tuning the coefficients in the linear span are done simultaneously during the neural network training process. Therefore, the neural network is given the opportunity to obtain a change of variable that transforms the coordinates so that they can be best fit to the (potentially) limited library provided by the user. This was on display with the R\"ossler, Kuramoto--Sivashinsky, and Mackey--Glass equations where our libraries only contained degree one and two monomials and the Gissinger system where we only use degree one, two, and three monomials. These examples are guided by intuition of what the mapping should be, but the process described above is best demonstrated on the Lorenz system. From years of previous work we have a good idea of what the Lorenz Poincar\'e mapping should look like in terms of symmetries and a jump discontinuity, but there is no indication that it should be (piecewise) polynomial. Our work shows that the network can transform the variables in such a way that they fit a piecewise quadratic map whose coefficients are tuned simultaneously. The result is a mapping guided by analytical rigour that is simple enough to gain specific insight from.}

Mappings that are conjugate to each other generate equivalent dynamics, thus forming an equivalence relation of the set of all continuous surjections of a topological space to itself. Hence, with the neural network we are able to classify and categorize distinct kinds of chaotic systems. For example, the R\"ossler, Kuramoto--Sivashinsky, and Mackey--Glass equations were all shown to be conjugate to the logistic map in certain parameter regions. Therefore, their chaotic attractors can be understood as a process of stretching and folding, as was made famous with Smale's horsehoe~\cite{Smale}. This lies in contrast to the chaotic switching between lobes observed in the Lorenz system. Beyond one-dimensional mappings lies a world of under-explored and little-understood chaotic dynamics. Therefore, the work initiated in this manuscript can help to classify the topology of chaotic attractors, while also providing simple mappings that can be analyzed with this goal in mind. In this way, we will be able to move toward better understanding not just chaotic, but hyperchaotic (more than one positive Lyapunov exponent) behaviour. \new{This also would include the study of periodically driven and Hamiltonian systems, which were notably absent from our discussion.} In the case of Hamiltonian systems, the discovered mappings must necessarily be measure-preserving. This property requires the targeted mapping to be at least two-dimensional and possess a specific structure, thus presenting an additional challenge to applying the method. This presents an obvious nontrivial extension of this method that we hope to report on in the future.     

Finally, the ability to obtain conjugacies of Poincar\'e mappings with actionable and interpretable nonlinear functions provides avenues for future computer-assisted analytical undertakings of chaotic systems. That is, a computer-assisted proof first showed that the Lorenz system is indeed chaotic using the Poincar\'e map~\cite{LorenzProof}. A more recent example of such an undertaking has coupled Poincar\'e mappings and computer-assisted proofs to show that the Kuramoto--Sivashinsky PDE is chaotic~\cite{KSProof}. Therefore, we conjecture that by using the neural network architecture of this manuscript one will be able to reach similar conclusions. Following similar computer-assisted proof methods, one could seek a desired threshold for the loss function of the neural network to prove that an exact conjugacy exists, followed by then proving that the simple conjugate mapping is chaotic. Since chaos is a topological invariant, it would then follow that the original system is also chaotic. If such a method of computer-assisted proof is possible, we suspect that it would considerably ease the difficulty of proving a system is chaotic.    
\section*{Acknowledgements}

J.J.B. thanks Joseph Bakarji for his help with learning TensorFlow 2.
SLB and JNK acknowledge funding support from the Air Force Office of Scientific Research (AFOSR FA9550-19-1-0386).  SLB acknowledges support from the Army Research Office (ARO W911NF-19-1-0045). JNK acknowledges support from the Air Force Office of Scientific Research (FA9550-19-1-0011).

\section*{Data Availability}
The data that support the findings of this study can be generated using the scripts in the repository {\bf github.com/jbramburger/Deep-Conjugacies}.

\setlength{\bibsep}{2.pt plus 1ex}
\begin{spacing}{.01}
	\small
\bibliographystyle{unsrt}
\bibliography{NN_conj}
\end{spacing}

\appendix

\section{Model Information}\label{app:1}

Here we report the specifications used to produce the models described throughout the manuscript. The first table collections the initial conditions used to generate the training data, as well as listing the location of the Poincar\'e section. The following table lists the network parameters found through parameter searches used to obtain the models. The 'Layers In' refers to the number of network layers to build $\hv$ and 'Layers Out' refers to the number of network layers to build $\hv^{-1}$. In every application these two values are always set to be the same, but this is not necessary.

\begin{center}
\footnotesize
\renewcommand{\arraystretch}{1.5}
\begin{tabular}{ |c||c|c|c|c|  }
 \hline
System Name & Reference & Initial Conditions & Section Location & Iterates \\
 \hline
 \hline
 R\"ossler   & \eqref{Rossler}    & $\begin{array}{l}
(0,-10,0),\ \ c = -9\\
(0,-15,0),\ \ \mathrm{otherwise}
\end{array}$ &  $\begin{array}{l}x_1 = 0\\\ \dot{x}_1 > 0\end{array}$ & $\begin{array}{l} 9877,\ \ c = 9 \\ 9845,\ \ c = 11 \\ 9788,\ \ c = 13 \\ 9738,\ \ c = 18 \\  \end{array}$  \\
  \hline
Lorenz &   \eqref{Lorenz}  & $(2,0,27)$   & $\begin{array}{l}x_3 = 27\\\ \dot{x}_3 < 0 \end{array}$ & 26631 \\
 \hline
 Gissinger & \eqref{Gissinger} & $(-1.5,1.5,1.3)$ & $\begin{array}{l}\hat x_1 + \hat x_2 = 0\\ \dot{\hat x}_1 + \dot{\hat x}_2 > 0\end{array}$ & 2615\\
  \hline
Kuramoto--Sivashinsky &   \eqref{Kuramoto}  & Random from $[0,0.1]^{14}$ & $\begin{array}{l}x_1 = 0\\\ \dot{x}_1 > 0\end{array}$ & 5728\\
 \hline
 Mackey--Glass & \eqref{Mackey}  & $x(t) = 0.5$, \ \ $t \in [-\tau,0]$   & $\begin{array}{l} (1)\ x(t) = x(t-\tau) \\ (2)\ \dot{x} = 0 \wedge \ddot{x} > 0 \end{array}$ & $\begin{array}{l} (1)\ 1739 \\ (2)\ 1739 \end{array}$\\
 \hline
\end{tabular}
\end{center}

\begin{center}
\footnotesize
\renewcommand{\arraystretch}{1.5}
\begin{tabular}{ |c||c|c|c|c|c|  }
 \hline
System Name & Layer Width & Layers In & Layers Out & Steps & Learning Rate \\
 \hline
 \hline
 R\"ossler ($c = 9$) & 100 & 1 & 1 & 2 & $5\times 10^{-3}$\\
 \hline
  R\"ossler ($c = 11$) & 80 & 1 & 1 & 2 & $5\times 10^{-3}$\\
 \hline
  R\"ossler ($c = 13$) & 80 & 1 & 1 & 2 & $5\times 10^{-3}$\\
\hline
  R\"ossler ($c = 18$) & 100 & 1 & 1 & 2 & $5\times 10^{-3}$\\
  \hline
  Lorenz & 200 & 3 & 3 & 2 & $1\times 10^{-3}$\\
  \hline
  Gissinger & 100 & 2 & 2 & 1 & $1\times 10^{-3}$\\
  \hline
 Kuramoto--Sivashinsky (1D) & 200 & 4 & 4 & 1 & $1\times 10^{-4}$\\
  \hline
 Kuramoto--Sivashinsky (2D) & 200 & 4 & 4 & 1 & $1\times 10^{-4}$\\
 \hline
 Mackey--Glass (1) & 200 & 4 & 4 & 2 & $5\times 10^{-4}$\\
 \hline
 Mackey--Glass (2) & 300 & 4 & 4 & 2 & $1\times 10^{-4}$\\
 \hline
\end{tabular}
\end{center}

\end{document}